\documentclass[11pt]{article}
\usepackage{amssymb}

\def\NZ{\mathbb{N}}
\def\GZ{\mathbb{Z}}

\def\OV#1{\overline {#1}}\def\op{\oplus}
\def\hang{\hangindent\parindent}
\def\textindent#1{\indent\llap{#1\enspace}\ignorespaces}
\def\re{\par\hang\textindent}

\def\QED{\hfill{$\Box$}}
\def \r{\rightarrow}

\def\mapdown#1{\llap{$\vcenter {\hbox {$\scriptstyle #1$}}$}
                                \Bigg\downarrow}
\def\mapdownr#1{\Bigg\downarrow\rlap{$\vcenter{\hbox
                                    {$\scriptstyle #1$}}$}}
\def\mapright#1#2{\smash{\mathop{\longrightarrow}\limits^{#1}_{#2}}}

\def\v5{\vskip .5truecm}

\def\HT{{\bf HT}}\def\LM{{\bf LM}}\def\LT{{\bf LT}}\def\KS{K\langle X\rangle}
\def\B{{\cal B}} \def\LC{{\bf LC}} \def\G{{\cal G}} \def\FRAC#1#2{\displaystyle{\frac{#1}{#2}}}
\def\SUM^#1_#2{\displaystyle{\sum^{#1}_{#2}}}  

\def\FRAC#1#2{\displaystyle{\frac{#1}{#2}}} 
\def\SUM^#1_#2{\displaystyle{\sum^{#1}_{#2}}}\def\G{{\cal G}}

\newdimen\DiagramCellHeight\DiagramCellHeight3em 
\newdimen\DiagramCellWidth\DiagramCellWidth3em 
\newdimen\MapBreadth\MapBreadth.04em 
\newdimen\MapShortFall\MapShortFall.4em 
\newdimen\PileSpacing\PileSpacing1em 
\def\labelstyle{\ifincommdiag\textstyle\else\scriptstyle\fi}
\let\objectstyle\displaystyle

\def\rTo{\HorizontalMap\empty-\empty-\rhvee}
\def\lTo{\HorizontalMap\lhvee-\empty-\empty}
\def\dTo{\VerticalMap\empty|\empty|\dhvee}
\def\uTo{\VerticalMap\uhvee|\empty|\empty}

\def\NW{\NorthWest\DiagonalMap{\lah111}{\laf100}{}{\laf100}{}(2,2)}
\def\NE{\NorthEast\DiagonalMap{\lah22}{\laf0}{}{\laf0}{}(2,2)}
\def\SW{\SouthWest\DiagonalMap{}{\laf0}{}{\laf0}{\lah11}(2,2)}
\def\SE{\SouthEast\DiagonalMap{}{\laf100}{}{\laf100}{\lah122}(2,2)}

\def\SSE{\SouthEast\DiagonalMap{}{\laf101}{}{\laf101}{\lah125}(2,4)}

\def\ESE{\SouthEast\DiagonalMap{}{\laf110}{}{\laf110}{\lah152}(4,2)}

\font\tenln=line10

\mathchardef\lt="313C \mathchardef\gt="313E

\def\rhvee{\mkern-10mu\gt}
\def\lhvee{\lt\mkern-10mu}
\def\dhvee{\vbox\tozpt{\vss\hbox{$\vee$}\kern0pt}}
\def\uhvee{\vbox\tozpt{\hbox{$\wedge$}\vss}}
\def\dhcvee{\vbox\tozpt{\vss\hbox{$\curlyvee$}\kern0pt}}
\def\dhvvee{\vbox\tozpt{\vss\hbox{$\vee$}\kern-.6ex\hbox{$\vee$}\kern0pt}}
\def\uhvvee{\vbox\tozpt{\hbox{$\wedge$}\kern-.6ex\hbox{$\wedge$}\vss}}
\def\twoheaddownarrow{\rlap{$\downarrow$}\raise-.5ex\hbox{$\downarrow$}}
\def\twoheaduparrow{\rlap{$\uparrow$}\raise.5ex\hbox{$\uparrow$}}
\def\rhla{\vbox\tozpt{\vss\hbox\tozpt{\hss\tenln\char'55}\kern\axisheight}}
\def\lhla{\vbox\tozpt{\vss\hbox\tozpt{\tenln\char'33\hss}\kern\axisheight}}
\def\rthooka{\raise.603ex\hbox{$\scriptscriptstyle\subset$}}
\def\lthooka{\raise.603ex\hbox{$\scriptscriptstyle\supset$}}
\def\rthookb{\raise-.022ex\hbox{$\scriptscriptstyle\subset$}}
\def\lthookb{\raise-.022ex\hbox{$\scriptscriptstyle\supset$}}

\def\SEpbk{\rlap{\smash{\kern0.1em \vrule depth 2.67ex height -2.55ex width 0%
.9em \vrule height -0.46ex depth 2.67ex width .05em }}}
\def\SWpbk{\llap{\smash{\vrule height -0.46ex depth 2.67ex width .05em \vrule
depth 2.67ex height -2.55ex width .9em \kern0.1em }}}
\def\NEpbk{\rlap{\smash{\kern0.1em \vrule depth -3.48ex height 3.67ex width 0%
.95em \vrule height 3.67ex depth -1.39ex width .05em }}}
\def\NWpbk{\llap{\smash{\vrule height 3.6ex depth -1.39ex width .05em \vrule
depth -3.48ex height 3.67ex width .95em \kern0.1em }}}

\newcount\cdna\newcount\cdnb\newcount\cdnc\newcount\cdnd\cdna=\catcode`\@%
\catcode`\@=11 \let\then\relax\def\loopa#1\repeat{\def\bodya{#1}\iteratea}%
\def\iteratea{\bodya\let\next\iteratea\else\let\next\relax\fi\next}\def\loopb
#1\repeat{\def\bodyb{#1}\iterateb}\def\iterateb{\bodyb\let\next\iterateb\else
\let\next\relax\fi\next} \def\mapctxterr{\message{commutative diagram: map
context error}}\def\mapclasherr{\message{commutative diagram: clashing maps}}%
\def\ObsDim#1{\expandafter\message{! diagrams Warning: Dimension \string#1 is
obsolete (ignored)}\global\let#1\ObsDimq\ObsDimq}\def\ObsDimq{\dimen@=}\def
\HorizontalMapLength{\ObsDim\HorizontalMapLength}\def\VerticalMapHeight{%
\ObsDim\VerticalMapHeight}\def\VerticalMapDepth{\ObsDim\VerticalMapDepth}\def
\VerticalMapExtraHeight{\ObsDim\VerticalMapExtraHeight}\def
\VerticalMapExtraDepth{\ObsDim\VerticalMapExtraDepth}\def\ObsCount#1{%
\expandafter\message{! diagrams Warning: Count \string#1 is obsolete (ignored%
)}\global\let#1\ObsCountq\ObsCountq}\def\ObsCountq{\count@=}\def
\DiagonalLineSegments{\ObsCount\DiagonalLineSegments}\def\tozpt{to\z@}\def
\sethorizhtdp{\dimen8=\axisheight\dimen9=\MapBreadth\advance\dimen8.5\dimen9%
\advance\dimen9-\dimen8}\def\horizhtdp{height\dimen8 depth\dimen9 }\def
\axisheight{\fontdimen22\the\textfont2 }\countdef\boxc@unt=14

\def\bombparameters{\hsize\z@\rightskip\z@ plus1fil minus\maxdimen
\parfillskip\z@\linepenalty9000 \looseness0 \hfuzz\maxdimen\hbadness10000
\clubpenalty0 \widowpenalty0 \displaywidowpenalty0 \interlinepenalty0
\predisplaypenalty0 \postdisplaypenalty0 \interdisplaylinepenalty0
\interfootnotelinepenalty0 \floatingpenalty0 \brokenpenalty0 \everypar{}%
\leftskip\z@\parskip\z@\parindent\z@\pretolerance10000 \tolerance10000
\hyphenpenalty10000 \exhyphenpenalty10000 \binoppenalty10000 \relpenalty10000
\adjdemerits0 \doublehyphendemerits0 \finalhyphendemerits0 \prevdepth\z@}\def
\startbombverticallist{\hbox{}\penalty1\nointerlineskip}

\def\pushh#1\to#2{\setbox#2=\hbox{\box#1\unhbox#2}}\def\pusht#1\to#2{\setbox#%
2=\hbox{\unhbox#2\box#1}}

\newif\ifallowhorizmap\allowhorizmaptrue\newif\ifallowvertmap
\allowvertmapfalse\newif\ifincommdiag\incommdiagfalse

\def\diagram{\hbox\bgroup$\vcenter\bgroup\startbombverticallist
\incommdiagtrue\baselineskip\DiagramCellHeight\lineskip\z@\lineskiplimit\z@
\mathsurround\z@\tabskip\z@\let\\\diagcr\allowhorizmaptrue\allowvertmaptrue
\halign\bgroup\lcdtempl##\rcdtempl&&\lcdtempl##\rcdtempl\cr}\def\enddiagram{%
\crcr\egroup\reformatmatrix\egroup$\egroup}

\def\lcdtempl{\futurelet\thefirsttoken\dolcdtempl}\newif\ifemptycell\def
\dolcdtempl{\ifx\thefirsttoken\rcdtempl\then\hskip1sp plus 1fil \emptycelltrue
\else\hfil$\emptycellfalse\objectstyle\fi}\def\rcdtempl{\ifemptycell\else$%
\hfil\fi}\def\diagcr{\cr} \def\across#1{\span\omit\mscount=#1 \loop\ifnum
\mscount>2 \spAn\repeat\ignorespaces}\def\spAn{\relax\span\omit\advance \mscount by
-1}

\def\CellSize{\afterassignment\cdhttowd\DiagramCellHeight}\def\cdhttowd{%
\DiagramCellWidth\DiagramCellHeight}\def\MapsAbut{\MapShortFall\z@}

\newcount\cdvdl\newcount\cdvdr\newcount\cdvd\newcount\cdbfb\newcount\cdbfr
\newcount\cdbfl\newcount\cdvdr\newcount\cdvdl\newcount\cdvd

\def\reformatmatrix{\bombparameters\cdvdl=\insc@unt\cdvdr=\cdvdl\cdbfb=%
\boxc@unt\advance\cdbfb1 \cdbfr=\cdbfb\setbox1=\vbox{}\dimen2=\z@\loop\setbox
0=\lastbox\ifhbox0 \dimen1=\lastskip\unskip\dimen5=\ht0 \advance\dimen5 \dimen
1 \dimen4=\dp0 \penalty1 \reformatrow\unpenalty\ht4=\dimen5 \dp4=\dimen4 \ht3%
\z@\dp3\z@\setbox1=\vbox{\box4 \nointerlineskip\box3 \nointerlineskip\unvbox1%
}\dimen2=\dimen1 \repeat\unvbox1}

\newif\ifcontinuerow

\def\reformatrow{\cdbfl=\cdbfr\noindent\unhbox0 \loopa\unskip\setbox\cdbfl=%
\lastbox\ifhbox\cdbfl\advance\cdbfl1\repeat\par\unskip\dimen6=2%
\DiagramCellWidth\dimen7=-\DiagramCellWidth\setbox3=\hbox{}\setbox4=\hbox{}%
\setbox7=\box\voidb@x\cdvd=\cdvdl\continuerowtrue\loopa\advance\cdvd-1
\adjustcells\ifcontinuerow\advance\dimen6\wd\cdbfl\cdda=.5\dimen6 \ifdim\cdda
<\DiagramCellWidth\then\dimen6\DiagramCellWidth\advance\dimen6-\cdda
\nopendvert\cdda\DiagramCellWidth\fi\advance\dimen7\cdda\dimen6=\wd\cdbfl
\reformatcell\advance\cdbfl-1 \repeat\advance\dimen7.5\dimen6 \outHarrow} \def
\adjustcells{\ifnum\cdbfr>\cdbfl\then\ifnum\cdvdr>\cdvd\then\continuerowfalse
\else\setbox\cdbfl=\hbox to\wd\cdvd{\lcdtempl\VonH{}\rcdtempl}\fi\else\ifnum
\cdvdr>\cdvd\then\advance\cdvdr-1 \setbox\cdvd=\vbox{}\wd\cdvd=\wd\cdbfl\dp
\cdvd=\dp1 \fi\fi}

\def\reformatcell{\sethorizhtdp\noindent\unhbox\cdbfl\skip0=\lastskip\unskip
\par\ifcase\prevgraf\reformatempty\or\reformatobject\else\reformatcomplex\fi
\unskip}\def\reformatobject{\setbox6=\lastbox\unskip\vadjdon6\outVarrow
\setbox6=\hbox{\unhbox6}\advance\dimen7-.5\wd6 \outHarrow\dimen7=-.5\wd6
\pusht6\to4}\newcount\globnum

\def\reformatcomplex{\setbox6=\lastbox\unskip\setbox9=\lastbox\unskip\setbox9%
=\hbox{\unhbox9 \skip0=\lastskip\unskip\global\globnum\lastpenalty\hskip\skip 0
}\advance\globnum9999 \ifcase\globnum\reformathoriz\or\reformatpile\or
\reformatHonV\or\reformatVonH\or\reformatvert\or\reformatHmeetV\fi}

\def\reformatempty{\vpassdon\ifdim\skip0>\z@\then\hpassdon\else\ifvoid2 \then
\else\advance\dimen7-.5\dimen0 \cdda=\wd2\advance\cdda.5\dimen0\wd2=\cdda\fi
\fi}\def\VonH{\doVonH6}\def\doVonH#1{\cdna-999#1\futurelet\thenexttoken
\dooVonH}\def\dooVonH{\let\next\relax\sethorizhtdp\ifallowhorizmap
\ifallowvertmap\then\ifx\thenexttoken[\then\let\next\VonHstrut\else
\sethorizhtdp\dimen0\MapBreadth\let\next\VonHnostrut\fi\else\mapctxterr\fi
\else\mapctxterr\fi\next}\def\VonHstrut[#1]{\setbox0=\hbox{$#1$}\dimen0\wd0%
\dimen8\ht0\dimen9\dp0 \VonHnostrut}\def\VonHnostrut{\setbox0=\hbox{}\ht0=%
\dimen8\dp0=\dimen9\wd0=.5\dimen0 \copy0\penalty\cdna\box0 \allowhorizmapfalse
\allowvertmapfalse}\def\reformatHonV{\hpassdon\doreformatHonV}\def
\reformatHmeetV{\dimen@=\wd9 \advance\dimen7-\wd9 \outHarrow\setbox6=\hbox{%
\unhbox6}\dimen7-\wd6 \advance\dimen@\wd6 \setbox6=\hbox to\dimen@{\hss}%
\pusht6\to4\doreformatHonV}\def\doreformatHonV{\setbox9=\hbox{\unhbox9 \unskip
\unpenalty\global\setbox\globbox=\lastbox}\vadjdon\globbox\outVarrow}\def
\reformatVonH{\vpassdon\advance\dimen7-\wd9 \outHarrow\setbox6=\hbox{\unhbox6%
}\dimen7=-\wd6 \setbox6=\hbox{\kern\wd9 \kern\wd6}\pusht6\to4}\def\hpassdon{}%
\def\vpassdon{\dimen@=\dp\cdvd\advance\dimen@\dimen4 \advance\dimen@\dimen5
\dp\cdvd=\dimen@\nopendvert}\def\vadjdon#1{\dimen8=\ht#1 \dimen9=\dp#1 }

\def\HorizontalMap#1#2#3#4#5{\sethorizhtdp\setbox1=\makeharrowpart{#1}\def
\arrowfillera{#2}\def\arrowfillerb{#4}\setbox5=\makeharrowpart{#5}\ifx
\arrowfillera\justhorizline\then\def\arra{\hrule\horizhtdp}\def\kea{\kern-0.%
01em}\let\arrstruthtdp\horizhtdp\else\def\kea{\kern-0.15em}\setbox2=\hbox{%
\kea${\arrowfillera}$\kea}\def\arra{\copy2}\def\arrstruthtdp{height\ht2 depth%
\dp2 }\fi\ifx\arrowfillerb\justhorizline\then\def\arrb{\hrule\horizhtdp}\def
\keb{kern-0.01em}\ifx\arrowfillera\empty\then\let\arrstruthtdp\horizhtdp\fi
\else\def\keb{\kern-0.15em}\setbox4=\hbox{\keb${\arrowfillerb}$\keb}\def\arrb
{\copy4}\ifx\arrowfilera\empty\then\def\arrstruthtdp{height\ht4 depth\dp4 }%
\fi\fi\setbox3=\makeharrowpart{{#3}\vrule width\z@\arrstruthtdp}%
\ifallowhorizmap\then\let\execmap\execHorizontalMap\else\let\execmap
\mapctxterr\fi\allowhorizmapfalse\gettwoargs}\def\makeharrowpart#1{\hbox{%
\mathsurround\z@\edef\next{#1}\ifx\next\empty\else$\mkern-1.5mu{\next}\mkern-%
1.5mu$\fi}}\def\justhorizline{-}

\def\execHorizontalMap{\dimen0=\wd6 \ifdim\dimen0<\wd7\then\dimen0=\wd7\fi
\dimen3=\wd3 \ifdim\dimen0<2em\then\dimen0=2em\fi\skip2=.5\dimen0 \ifincommdiag
plus 1fill\fi minus\z@\advance\skip2-.5\dimen3 \skip4=\skip2 \advance\skip2-\wd1
\advance\skip4-\wd5 \kern\MapShortFall\box1 \xleaders
\arra\hskip\skip2 \vbox{\lineskiplimit\maxdimen\lineskip.5ex \ifhbox6 \hbox to%
\dimen3 {\hss\box6\hss}\fi\vtop{\box3 \ifhbox7 \hbox to\dimen3 {\hss\box7\hss
}\fi}}\ifincommdiag\kern-.5\dimen3\penalty-9999\null\kern.5\dimen3\fi
\xleaders\arrb\hskip\skip4 \box5 \kern\MapShortFall}

\def\reformathoriz{\vadjdon6\outVarrow\ifvoid7\else\mapclasherr\fi\setbox2=%
\box9 \wd2=\dimen7 \dimen7=\z@\setbox7=\box6 }

\def\resetharrowpart#1#2{\ifvoid#1\then\ifdim#2=\z@\else\setbox4=\hbox{%
\unhbox4\kern#2}\fi\else\ifhbox#1\then\setbox#1=\hbox to#2{\unhbox#1}\else
\widenpile#1\fi\pusht#1\to4\fi}\def\outHarrow{\resetharrowpart2{\wd2}\pusht2%
\to4\resetharrowpart7{\dimen7}\pusht7\to4\dimen7=\z@}

\def\pile#1{{\incommdiagtrue\let\pile\innerpile\allowvertmapfalse
\allowhorizmaptrue\baselineskip.5\PileSpacing\lineskip\z@\lineskiplimit\z@
\mathsurround\z@\tabskip\z@\let\\\pilecr\vcenter{\halign{\hfil$##$\hfil\cr#1
\crcr}}}\ifincommdiag\then\ifallowhorizmap\then\penalty-9998
\allowvertmapfalse\allowhorizmapfalse\else\mapctxterr\fi\fi}\def\pilecr{\cr}%
\def\innerpile#1{\noalign{\halign{\hfil$##$\hfil\cr#1 \crcr}}}

\def\reformatpile{\vadjdon9\outVarrow\ifvoid7\else\mapclasherr\fi\penalty1
\setbox9=\hbox{\unhbox9 \unskip\unpenalty\setbox9=\lastbox\unhbox9 \global
\setbox\globbox=\lastbox}\unvbox\globbox\setbox9=\vbox{}\setbox7=\vbox{}%
\loopb\setbox6=\lastbox\ifhbox6 \skip3=\lastskip\unskip\splitpilerow\repeat
\unpenalty\setbox9=\hbox{$\vcenter{\unvbox9}$}\setbox2=\box9 \dimen7=\z@}\def
\pilestrut{\vrule height\dimen0 depth\dimen3 width\z@}\def\splitpilerow{%
\dimen0=\ht6 \dimen3=\dp6 \noindent\unhbox6\unskip\setbox6=\lastbox\unskip
\unhbox6\par\setbox6=\lastbox\unskip\ifcase\prevgraf\or\setbox6=\hbox\tozpt{%
\hss\unhbox6\hss}\ht6=\dimen0 \dp6=\dimen3 \setbox9=\vbox{\vskip\skip3 \hbox
to\dimen7{\hfil\box6}\nointerlineskip\unvbox9}\setbox7=\vbox{\vskip\skip3
\hbox{\pilestrut\hfil}\nointerlineskip\unvbox7}\or\setbox7=\vbox{\vskip\skip3
\hbox{\pilestrut\unhbox6}\nointerlineskip\unvbox7}\setbox6=\lastbox\unskip
\setbox9=\vbox{\vskip\skip3 \hbox to\dimen7{\pilestrut\unhbox6}%
\nointerlineskip\unvbox9}\fi\unskip}

\def\widenpile#1{\setbox#1=\hbox{$\vcenter{\unvbox#1 \setbox8=\vbox{}\loopb
\setbox9=\lastbox\ifhbox9 \skip3=\lastskip\unskip\setbox8=\vbox{\vskip\skip3 \hbox
to\dimen7{\unhbox9}\nointerlineskip\unvbox8}\repeat\unvbox8 }$}}

\def\justverticalline{|}\def\makevarrowpart#1{\hbox to\MapBreadth{\hss$\kern
\MapBreadth{#1}$\hss}}\def\VerticalMap#1#2#3#4#5{\setbox1=\makevarrowpart{#1}%
\def\arrowfillera{#2}\setbox3=\makevarrowpart{#3}\def\arrowfillerb{#4}\setbox
5=\makevarrowpart{#5}\ifx\arrowfillera\justverticalline\then\def\arra{\vrule
width\MapBreadth}\def\kea{\kern-0.05ex}\else\def\kea{\kern-0.35ex}\setbox2=%
\vbox{\kea\makevarrowpart\arrowfillera\kea}\def\arra{\copy2}\fi\ifx
\arrowfillerb\justverticalline\then\def\arrb{\vrule width\MapBreadth}\def\keb
{\kern-0.05ex}\else\def\keb{\kern-0.35ex}\setbox4=\vbox{\keb\makevarrowpart
\arrowfillerb\keb}\def\arrb{\copy4}\fi\ifallowvertmap\then\let\execmap
\execVerticalMap\else\let\execmap\mapctxterr\fi\allowhorizmapfalse\gettwoargs }

\def\execVerticalMap{\setbox3=\makevarrowpart{\box3}\setbox0=\hbox{}\ht0=\ht3
\dp0\z@\ht3\z@\box6 \setbox8=\vtop spread2ex{\offinterlineskip\box3 \xleaders
\arrb\vfill\box5 \kern\MapShortFall}\dp8=\z@\box8 \kern-\MapBreadth\setbox8=%
\vbox spread2ex{\offinterlineskip\kern\MapShortFall\box1 \xleaders\arra\vfill
\box0}\ht8=\z@\box8 \ifincommdiag\then\kern-.5\MapBreadth\penalty-9995 \null
\kern.5\MapBreadth\fi\box7\hfil}

\newcount\colno\newdimen\cdda\newbox\globbox\def\reformatvert{\setbox6=\hbox{%
\unhbox6}\cdda=\wd6 \dimen3=\dp\cdvd\advance\dimen3\dimen4 \setbox\cdvd=\hbox
{}\colno=\prevgraf\advance\colno-2 \loopb\setbox9=\hbox{\unhbox9 \unskip
\unpenalty\dimen7=\lastkern\unkern\global\setbox\globbox=\lastbox\advance
\dimen7\wd\globbox\advance\dimen7\lastkern\unkern\setbox9=\lastbox\vtop to%
\dimen3{\unvbox9}\kern\dimen7 }\ifnum\colno>0 \ifdim\wd9<\PileSpacing\then
\setbox9=\hbox to\PileSpacing{\unhbox9}\fi\dimen0=\wd9 \advance\dimen0-\wd
\globbox\setbox\cdvd=\hbox{\kern\dimen0 \box\globbox\unhbox\cdvd}\pushh9\to6%
\advance\colno-1 \setbox9=\lastbox\unskip\repeat\advance\dimen7-.5\wd6
\advance\dimen7.5\cdda\advance\dimen7-\wd9 \outHarrow\dimen7=-.5\wd6 \advance
\dimen7-.5\cdda\pusht9\to4\pusht6\to4\nopendvert\dimen@=\dimen6\advance
\dimen@-\wd\cdvd\advance\dimen@-\wd\globbox\divide\dimen@2 \setbox\cdvd=\hbox
{\kern\dimen@\box\globbox\unhbox\cdvd\kern\dimen@}\dimen8=\dp\cdvd\advance
\dimen8\dimen5 \dp\cdvd=\dimen8 \ht\cdvd=\z@}

\def\outVarrow{\ifhbox\cdvd\then\deepenbox\cdvd\pusht\cdvd\to3\else
\nopendvert\fi\dimen3=\dimen5 \advance\dimen3-\dimen8 \setbox\cdvd=\vbox{%
\vfil}\dp\cdvd=\dimen3} \def\nopendvert{\setbox3=\hbox{\unhbox3\kern\dimen6}}%
\def\deepenbox\cdvd{\setbox\cdvd=\hbox{\dimen3=\dimen4 \advance\dimen3-\dimen
9 \setbox6=\hbox{}\ht6=\dimen3 \dp6=-\dimen3 \dimen0=\dp\cdvd\advance\dimen0%
\dimen3 \unhbox\cdvd\dimen3=\lastkern\unkern\setbox8=\hbox{\kern\dimen3}%
\loopb\setbox9=\lastbox\ifvbox9 \setbox9=\vtop to\dimen0{\copy6
\nointerlineskip\unvbox9 }\dimen3=\lastkern\unkern\setbox8=\hbox{\kern\dimen3%
\box9\unhbox8}\repeat\unhbox8 }}

\newif\ifPositiveGradient\PositiveGradienttrue\newif\ifClimbing\Climbingtrue
\newcount\DiagonalChoice\DiagonalChoice1 \newcount\lineno\newcount\rowno
\newcount\charno\def\laf{\afterassignment\xlaf\charno='}\def\xlaf{\hbox{%
\tenln\char\charno}}\def\lah{\afterassignment\xlah\charno='}\def\xlah{\hbox{%
\tenln\char\charno}}\def\makedarrowpart#1{\hbox{\mathsurround\z@${#1}$}}\def
\lad{\afterassignment\xlad\charno='}\def\xlad{\setbox2=\xlaf\setbox0=\hbox to%
.5\wd2{$\hss\ldot\hss$}\ht0=.25\ht2 \dp0=\ht0 \hbox{\mv-\ht0\copy0 \mv\ht0%
\box0}}

\def\DiagonalMap#1#2#3#4#5{\ifPositiveGradient\then\let\mv\raise\else\let\mv
\lower\fi\setbox2=\makedarrowpart{#2}\setbox1=\makedarrowpart{#1}\setbox4=%
\makedarrowpart{#4}\setbox5=\makedarrowpart{#5}\setbox3=\makedarrowpart{#3}%
\let\execmap\execDiagonalLine\gettwoargs}

\def\makeline#1(#2,#3;#4){\hbox{\dimen1=#2\relax\dimen2=#3\relax\dimen5=#4%
\relax\vrule height\dimen5 depth\z@ width\z@\setbox8=\hbox to\dimen1{\tenln#1%
\hss}\cdna=\dimen5 \divide\cdna\dimen2 \ifnum\cdna=0 \then\box8 \else\dimen4=%
\dimen5 \advance\dimen4-\dimen2 \divide\dimen4\cdna\dimen3=\dimen1 \cdnb=%
\dimen2 \divide\cdnb1000 \divide\dimen3\cdnb\cdnb=\dimen4 \divide\cdnb1000
\multiply\dimen3\cdnb\dimen6\dimen1 \advance\dimen6-\dimen3 \cdnb=0
\ifPositiveGradient\then\dimen7\z@\else\dimen7\cdna\dimen4 \multiply\dimen4-1
\fi\loop\raise\dimen7\copy8 \ifnum\cdnb<\cdna\hskip-\dimen6 \advance\cdnb1
\advance\dimen7\dimen4 \repeat\fi}}\newdimen\objectheight\objectheight1.5ex

\def\execDiagonalLine{\setbox0=\hbox\tozpt{\cdna=\xcoord\cdnb=\ycoord\dimen8=%
\wd2 \dimen9=\ht2 \dimen0=\cdnb\DiagramCellHeight\advance\dimen0-2%
\MapShortFall\advance\dimen0-\objectheight\setbox2=\makeline\box2(\dimen8,%
\dimen9;.5\dimen0)\setbox4=\makeline\box4(\dimen8,\dimen9;.5\dimen0)\dimen0=2%
\wd2 \advance\dimen0-\cdna\DiagramCellWidth\advance\dimen0 2\DiagramCellWidth
\dimen2\DiagramCellHeight\advance\dimen2-\MapShortFall\dimen1\dimen2 \advance
\dimen1-\ht1 \advance\dimen2-\ht2 \dimen6=\dimen2 \advance\dimen6.25\dimen8
\dimen3\dimen2 \advance\dimen3-\ht3 \dimen4=\dimen2 \dimen7=\dimen2 \advance
\dimen4-\ht4 \advance\dimen7-\ht7 \advance\dimen7-.25\dimen8
\ifPositiveGradient\then\hss\raise\dimen4\hbox{\rlap{\box5}\box4}\llap{\raise
\dimen6\box6\kern.25\dimen9}\else\kern-.5\dimen0 \rlap{\raise\dimen1\box1}%
\raise\dimen2\box2 \llap{\raise\dimen7\box7\kern.25\dimen9}\fi\raise\dimen3%
\hbox\tozpt{\hss\box3\hss}\ifPositiveGradient\then\rlap{\kern.25\dimen9\raise
\dimen7\box7}\raise\dimen2\box2\llap{\raise\dimen1\box1}\kern-.5\dimen0 \else
\rlap{\kern.25\dimen9\raise\dimen6\box6}\raise\dimen4\hbox{\box4\llap{\box5}}%
\hss\fi}\ht0\z@\dp0\z@\box0}

\def\NorthWest{\PositiveGradientfalse\Climbingtrue\DiagonalChoice0 }\def
\NorthEast{\PositiveGradienttrue\Climbingtrue\DiagonalChoice1 }\def\SouthWest
{\PositiveGradienttrue\Climbingfalse\DiagonalChoice3 }\def\SouthEast{%
\PositiveGradientfalse\Climbingfalse\DiagonalChoice2 }

\newif\ifmoremapargs\def\gettwoargs{\setbox7=\box\voidb@x\setbox6=\box
\voidb@x\moremapargstrue\def\whichlabel{6}\def\xcoord{2}\def\ycoord{2}\def
\contgetarg{\def\whichlabel{7}\ifmoremapargs\then\let\next\getanarg\let
\contgetarg\execmap\else\let\next\execmap\fi\next}\getanarg}\def\getanarg{%
\futurelet\thenexttoken\switcharg}\def\getlabel#1#2#3{\setbox#1=\hbox{$%
\labelstyle\>{#3}\>$}\dimen0=\ht#1\advance\dimen0 .4ex\ht#1=\dimen0 \dimen0=%
\dp#1\advance\dimen0 .4ex\dp#1=\dimen0 \contgetarg}\def\eatspacerepeat{%
\afterassignment\getanarg\let\junk= }\def\catcase#1:{{\ifcat\noexpand
\thenexttoken#1\then\global\let\xcase\docase\fi}\xcase}\def\tokcase#1:{{\ifx
\thenexttoken#1\then\global\let\xcase\docase\fi}\xcase}\def\default:{\docase}%
\def\docase#1\esac#2\esacs{#1}\def\skipcase#1\esac{}\def\getcoordsrepeat(#1,#%
2){\def\xcoord{#1}\def\ycoord{#2}\getanarg}\let\esacs\relax\def\switcharg{%
\global\let\xcase\skipcase\catcase{&}:\moremapargsfalse\contgetarg\esac
\catcase\bgroup:\getlabel\whichlabel-\esac\catcase^:\getlabel6\esac\catcase_:%
\getlabel7\esac\tokcase{~}:\getlabel3\esac\tokcase(:\getcoordsrepeat\esac \catcase{
}:\eatspacerepeat\esac\default:\moremapargsfalse\contgetarg\esac \esacs}

\catcode`\@=\cdna

\title{Algebras and their Associated Monomial Algebras\thanks{Project supported by
the National Natural Science Foundation of China (10571038).
\newline e-mail: huishipp@yahoo.com}} \vskip 1truecm
\author{Huishi Li \\
{\small Department of Applied Mathematics}\\
{\small College of Information Science and Technology}\\
{\small Hainan University}\\
{\small  Haikou 570228, China}}

\date{}

\begin{document}
\maketitle
\begin{center}
\begin{minipage}{120mm}
{\small {\bf Abstract.} Let $R=\oplus_{\Gamma\in\Gamma}R_{\gamma}$ be a
$\Gamma$-graded $K$-algebra over a field $K$, where $\Gamma$ is a totally ordered
semigroup, and let $I$ be an ideal of $R$. Considering the $\Gamma$-grading
filtration $FR$ of $R$ and the $\Gamma$-filtration $FA$ induced by $FR$ for the
quotient $K$-algebra $A=R/I$, we show that there is a $\Gamma$-graded $K$-algebra
isomorphism $G(A)\cong \overline{ A}=R/\langle\hbox{\bf HT} (I)\rangle$, where
$G(A)$ is the associated $\Gamma$-graded $K$-algebra of $A$ defined by $FA$, and
$\langle\hbox{\bf HT}(I)\rangle$ is the $\Gamma$-graded ideal of $R$ generated by
the set of head terms of $I$. In the case that $\Gamma$ is an ordered monoid with a
well-ordering, this result enables us to lift many nice structural  properties of
$\overline{A}$ to $A$ theoretically, and  the natural connection with Gr\"obner
basis theory leads to effective realization lifting information from the associated
monomial algebras in both commutative and noncommutative cases.}
\end{minipage}\end{center}{\parindent=0pt\vskip 6pt
{\bf 2000 Mathematics Classification:} Primary 16W70; Secondary 16Z05.\vskip 6pt
{\bf Key words} Filtration, gradation, monomial algebra, Gr\"obner basis}

\vskip .5truecm
\section*{0. Introduction}
Let $K$ be a field, and let $R$ be one of the following
$K$-algebras:{\parindent=.5truecm\vskip 6pt \re{$\bullet$} $K[x_1,...,x_n]$, the
commutative polynomial $K$-algebra in $n$ variables, which has the standard
$K$-basis $\B =\{ x_1^{\alpha_1}\cdots x_n^{\alpha_n}~|~ \alpha_i\in\NZ\}$;
\re{$\bullet$} $K[a_1,...,a_n]$, the $K$-algebra generated by $a_1,...,a_n$ subject
to the relations
$$a_ja_i=\lambda_{ji}a_ia_j,\quad \lambda_{ji}\in K^*,~1\le i<j\le n,$$
which has the standard $K$-basis $\B=\{ a_1^{\alpha_1}\cdots
a_n^{\alpha_n}~|~\alpha_i\in\NZ\}$. (It is known that in computational algebra this
algebra is studied as a typical solvable polynomial algebra (e.g., [K-RW], [Li1]),
and in the study of noncommutative algebraic geometry it is called the coordinate
ring of the $n$-dimensional quantum affine $K$-space.);
\re{$\bullet$} $K\langle X_1,...,X_n\rangle$, the noncommutative free $K$-algebra
generated by $X=\{ X_1,...,X_n\}$, which has the standard $K$-basis $\B =\{
1,~X_{j_1}\cdots X_{j_s}~|~X_{j_i}\in X,~s\ge 1\}$;                \re{$\bullet$}
$KQ$, the path algebra defined by a finite directed graph $Q$ over $K$, which has
the standard $K$-basis $\B$ consisting of all finite paths including vertices as
paths of length 0. \re{~}}    {\parindent=0pt \par
Then, it is well-known that $R$ holds a well-developed (commutative or
noncommutative) Gr\"obner basis theory (cf. [Bu], [Mor], [K-RW], [Gr]). More
precisely, let $\prec$ be a monomial ordeirng on $\B$. If $I$ is an ideal of $R$,
then, there is a (finite or infinite) Gr\"obner basis $\G\subset I$ in the sense
that
$$\langle\LM (I)\rangle =\langle\LM (\G )\rangle ,$$
where $\langle\LM (I)\rangle$ is the  ideal of $R$ generated by the set of leading
monomials $\LM (I)$ of $I$ and $\langle\LM (\G )\rangle$ is the  ideal of $R$
generated by the set of  leading monomials $\LM (\G )$ of $\G$ (see section 5 for
the definition of a leading monomial). Put $A=R/I$ and $\OV A=R/\langle\LM
(I)\rangle$. In the literature the $K$-algebra $\OV A$ is usually called the {\it
associated monomial algebra} of the $K$-algebra $A$ due to the fact that
$\langle\LM (I)\rangle$ is a monomial ideal of $R$ (e.g., see  [An2], [{G-IL],
[G-I2], [GZ]). Historically, monomial algebras are studied and used widely in many
mathematical areas such as algebraic geometry, representation theory of algebras,
algebraic combinatorics,  as such algebras may be understood more easily, and
especially, may be manipulated on computer more effectively. To see the influence
of the monomial algebra $\OV A$ on the algebra $A$, a motive example, which can be
found in any computational work concerning Hilbert function, Hilbert series and
Poincar\'e series of a (graded) algebra, is worthwhile to be recalled here. Let $R$
be the free $K$-algebra $K\langle X_1,...,X_n\rangle$, and let $\prec$ be a
monomial ordering on $\B$. For an ideal $I$ of $R$, put $A=R/I$, $\OV
A=R/\langle\LM (I)\rangle$. Then the following statements
hold.}{\parindent=.75truecm\vskip 6pt                     \re{(1)} The image of the
set $\B -\LM (I)$ in $A$, respectively in $\OV A$, forms a $K$-basis for $A$,
respectively for $\OV A$.
\re{(2)}  With respect to the natural $\NZ$-filtration on $A$ and $\OV A$ (see
section 1 for the definition), $A$ and $\OV A$ have the same Hilbert function and
hence have the same growth, or equivalently, $A$ and $\OV A$ have the same
Gelfand-Kirillov dimension.
\re{(3)} If $I$ is an $\NZ$-graded ideal of $R$, then the $\NZ$-graded algebras $A$
and $\OV A$ have the same Hilbert series.
\re{(4)} If $\G$ is a Gr\"obner basis with respect to $(\B ,\prec)$, then all
invariants in (1) -- (3) are determined by $\LM (\G)\subset\langle\LM (I)\rangle$
and computable by means of some computer algebra system such as  BERGMAN
[CU].\parindent=0pt\re{~}  Similar results hold for other commonly studied algebras
that hold a Gr\"obner basis theory.}}\par
 Based on the above review, we are naturally concerned about the following problem.
\vskip 6pt {\parindent=0pt
{\bf Question} How to transfer as many as possible nice structural and
computational properties of $\OV A=R/\langle\LM (I)\rangle$ to $A=R/I$.}\vskip 6pt
In the case that $R$ is the free $K$-algebra $K\langle X_1,...,X_n\rangle$, if $I$
is an ideal of $R$, $A=R/I$ is the quotient algebra of $R$ defined by $I$, and
$G^{\NZ}(A)$ is the associated $\NZ$-graded $K$-algebra of $A$ with respect to its
natural $\NZ$-filtration $F^{\NZ}A$ induced by the $\NZ$-grading filtration
$F^{\NZ}R$ of $R$ defined by its natura $\NZ$-gradation (see section 1 for the
definition of a $\Gamma$-grading filtration), then it follows from ([Li1] Chapter
III Proposition 3.1) that $G^{\NZ}(A)\cong R/\langle\HT (I)\rangle$ as
$\mathbb{N}$-graded $K$-algebras, where $\langle\HT (I)\rangle$ is the
$\mathbb{N}$-graded ideal of $R$ generated by the set of head terms of $I$ (see
section 2 for the definition of a head term), and that if furthermore $\G$ is a
Gr\"obner basis for $I$ with respect to some graded monomial ordering on $\B$, then
$\langle\HT (I)\rangle =\langle\HT (\G )\rangle$. In [Li2], this result was
extended to propose a general PBW property for quotient algebras of a $\GZ$-graded
algebra, and for quotient algebras of a path algebra (including free algebra), a
solution to the general PBW problem is given by means of Gr\"obner bases.
Enlightened by [Li2], in the present paper we strive for a solution to the problem
posed above by virtue of the $\B$-filtration and Gr\"obner bases, but consideration
is made in a more general setting.  The contents of this paper are arranged as
follows.\vskip 6pt
1. $\Gamma$-filtered Algebras and Modules\par
2. With $\Gamma$-grading Filtration: $G(R/I)\cong R/\langle\HT (I)\rangle$\par
3. Basic Lifting Properties\par
4. Lifting Homological Properties\par
5. With Gr\"obner Bases: $G^{\B}(R/I)\cong  R/\langle\LM (\G )\rangle$ \&
$G^{\NZ}(R/I)\cong R/\langle\HT (\G )\rangle$\par
6. The First Application\par
7. Realization via Gr\"obner Bases and Ufnarovski Graphs{\parindent=0pt\vskip 6pt
{\bf Convention throughout the paper}\par
Let $K$ be a field. All algebras considered are associative $K$-algebras with
identity 1, and all modules, unless otherwise stated, are unitary left modules. Let
$R$ be a $K$-algebra and $S\subset R$. We write $\langle S\rangle$ for the
(two-sided) ideal of $R$ generated by the subset $S$, and write $\langle S]$ for
the left ideal of $R$ generated by $S$. Moreover, $K^*=K-\{ 0\}$.}\v5

Here we point out in advance that since the $\Gamma$-filtration is less studied in
the literature, and due to the nontrivial difference between a general ordered
semigroup $\Gamma$ and $\NZ$, we introduce this notion and the associated
$\Gamma$-graded structure in section 1 in a slightly detailed manner; besides,
although all results of sections 3 -- 4 are well-known in the case of $\Gamma
=\NZ$, to convince the reader, we provide a detailed proof for each result
concerning $\Gamma$-filtration, for, the author cannot say that all of them are
just a trivial imitation of the $\NZ$-filtered case.
 \v5
\section*{1. $\Gamma$-filtered Algebras and Modules}
In this section, $\Gamma$ denotes a {\it totally ordered} semigrouop, i.e.,
$\Gamma$ is a semigroup on which there is a total ordering $\prec$ that is
compatible with the binary operation of $\Gamma$ in the sense that for $\gamma_1$,
$\gamma_2$, $\gamma\in\Gamma$,
$$\gamma_1\prec\gamma_2~\hbox{implies}~\gamma\gamma_1\prec\gamma\gamma_2~\hbox{and}~
\gamma_1\gamma\prec\gamma_2\gamma .$$ {\parindent=0pt\par
{\bf $\Gamma$-filtered Algebra}\par
A $K$-algebra $A$ is said to be $\Gamma$-{\it filtered} if there is a family $FA=\{
F_{\gamma}A\}_{\gamma\in\Gamma}$ consisting of $K$-subspaces $F_{\gamma}A$ of $A$,
such that {\parindent=1.2truecm\vskip 6pt
\re{\bf (F1)} $A=\cup_{\gamma\in\Gamma}F_{\gamma}A$,\par
\re{\bf (F2)} $F_{\gamma_1}A\subseteq F_{\gamma_2}A$ if
$\gamma_1\preceq\gamma_2$,\par
\re{\bf (F3)} $F_{\gamma_1}AF_{\gamma_2}A\subseteq F_{\gamma_1\gamma_2}A$,
$\gamma_1,\gamma_2\in\Gamma$.}\vskip 6pt
If $\Gamma$ has a smallest element $\gamma_0$, we also ask that $1\in
F_{\gamma_0}A$.}\par
In the case that $\Gamma =\GZ$ and $A$ is a $\GZ$-filtered $K$-algebra, if $F_nA=\{
0\}$ for all $n<0$, then $A$ becomes an $\NZ$-filtered $K$-algebra, which is also
called a {\it positively filtered} $K$-algebra.\vskip 6pt
In the definition given above, the family $FA=\{ F_{\gamma}A\}_{\gamma\in\Gamma}$
is usually called a $\Gamma$-{\it filtration} of $A$. \v5
{\parindent=0pt\v5
{\bf Natural $\NZ$-filtration defined by lengths of monomials}\par
Let $A=K[T]$ be a $K$-algebra generated by $T=\{ a_i\}_{i\in J}$ over $K$. Then
each element $a\in A$ can be written as a finite sum of the form
$$a=\sum\lambda_{i_1\cdots i_s}a_{i_1}^{\alpha_1}\cdots a_{i_s}^{\alpha_s},\quad
a_{i_j}\in T,~\lambda_{i_1\cdots i_s}\in K,~\alpha_i,s\in\NZ,~s\ge 1.$$
By abusing language, a nonzero element of the form $u=a_{i_1}^{\alpha_1}\cdots
a_{i_s}^{\alpha_s}$ is called a {\it monomial} of $A$, and the {\it length} of
$u$, denoted $l(u)$, is defined as $l(u)=\alpha_1+\cdots +\alpha_s$. Let $\Omega$
be the set of all monomials in $A$, i.e., $\Omega =\{ u=a_{i_1}^{\alpha_1}\cdots
a_{i_s}^{\alpha_s}~|~a_{i_j}\in T,~\alpha_i,s\in\NZ,~s\ge 1\}$.  For each
$p\in\NZ$, let $F_pA$ denote the $K$-subspace of $A$ spanned by all monomials of
length less than or equal to $p$, that is
$$F_pA=K\hbox{-span}\left\{ u\in\Omega~\Big |~l(u)\le p\right\} .$$
It is easy to see that the family $FA=\{ F_pA\}_{p\in\NZ}$ satisfies the foregoing
conditions (F1)--(F3). This $\NZ$-filtration is called the {\it natural
$\NZ$-filtration} of $A$ defined by lengths of monomials.\v5
{\bf $\Gamma$-grading filtration}\par
Let  $R$ be a $\Gamma$-graded $K$-algebra, that is,
$R=\oplus_{\gamma\in\Gamma}R_{\gamma}$, where for each $\gamma\in\Gamma$,
$R_{\gamma}$ is a $K$-subspace of $R$, and for any $\gamma_1$, $\Gamma_2\in\Gamma$,
$R_{\gamma_1}R_{\gamma_2}\subseteq R_{\gamma_1\gamma_2}$.  Put
$$F_{\gamma}R=\bigoplus_{\gamma'\preceq\gamma}R_{\gamma'},\quad\gamma\in\Gamma .$$
Then it may be checked directly that the family $FR=\{
F_{\gamma}R\}_{\gamma\in\Gamma}$ satisfies the foregoing conditions (F1)--(F3).
This $\Gamma$-filtration is called the {\it $\Gamma$-grading filtration} of $R$
defined by the given $\Gamma$-gradation of $R$ \v5
{\bf Example} (1) Let Let $R$ be the free $K$-algebra $K\langle
X_1,...,X_n\rangle$, or the commutative polynomial $K$-algebra $K[x_1,...,x_n]$, or
the coordinate ring of the $n$-dimensional quantum affine $K$-space
$K[a_1,...,a_n]$, or the path algebra $KQ$ defined by a finite directed graph $Q$
over $K$, and let $\B$ be the standard $K$-basis of $R$.  Then $R$ is $\NZ$-graded
by the  $\NZ$-gradation $\{ R_p\}_{p\in\NZ}$ with $R_p=K\hbox{-span}\{
u\in\B~|~l(u)=p\}$. It is easy to see that the natural $\NZ$-filtration of $R$
defined by lengths of monomials coincides with the $\NZ$-grading filtration defined
by the $\NZ$-gradation $\{ R_p\}_{p\in\NZ}$.}
\par
If furthermore $\prec$ is a monomial ordering on $\B$, then $\B$ becomes an ordered
semigroup with the well-ordering $\prec$ (in the case that $R=KQ$, $\B\cup\{ 0\}$
is considered). It turns out that $R$ is $\B$-graded, i.e., $R=\oplus_{u\in\B}R_u$
with $R_u=Ku$, and consequently, this $\B$-gradation defines the $\B$-grading
filtration $FR=\{ F_uR\}_{u\in\B}$ of $R$ with $F_uR=\oplus_{u'\preceq
u}R_{u'}$.\par
Both the $\B$-filtration and the $\NZ$-filtration of $R$ will be used in later
section 5.{\parindent=0pt\v5
{\bf The associated $\Gamma$-graded algebra}\par
Let $A$ be a $\Gamma$-filtered $K$-algebra with  $\Gamma$-filtration $FA=\{
F_{\gamma}A\}_{\gamma\in\Gamma}$. Put
$$\begin{array}{l} F_{\gamma}^*A=\bigcup_{\gamma'\prec\gamma}F_{\gamma'}A,\quad \gamma\in\Gamma ,\\

\hbox{where}~F_{\gamma_0}^*A=\{ 0\}~\hbox{if}~A~\hbox{has a smallest element}~\gamma_0.\end{array}$$
The {\it associated $\Gamma$-graded $K$-algebra} of $A$, denoted $G(A)$, is defined
as
$$G(A)=\bigoplus_{\gamma\in\Gamma}G(A)_{\gamma}~\hbox{with}~G(A)_{\gamma}=F_{\gamma}A/F_{\gamma}^*A,$$
where the multiplication is defined by extending the maps
$$\begin{array}{ccc} G(A)_{\gamma_1}\times G(A)_{\gamma_2}&\mapright{}{}&G(A)_{\gamma_1\gamma_2}\\
(\OV{a_{\gamma_1}},\OV{a_{\gamma_2}})&\mapsto&\OV{a_{\gamma_1}a_{\gamma_2}}\end{array}$$
to $G(A)\times G(A)\r G(A)$, in which $\OV{a_{\gamma_1}}$, $\OV{a_{\gamma_2}}$ are
the images of $a_{\gamma_1}\in F_{\gamma_1}A$, $a_{\gamma_2}\in F_{\gamma_2}A$ in
$G(A)_{\gamma_1}=F_{\gamma_1}A/F_{\gamma_1}^*A$ and
$G(A)_{\gamma_2}=F_{\gamma_2}A/F_{\gamma_2}^*A$ respectively, and
$\OV{a_{\gamma_1}a_{\gamma_2}}$ is the image of $a_{\gamma_1}a_{\gamma_2}\in
F_{\gamma_1\gamma_2}A$ in
$G(A)_{\gamma_1\gamma_2}=F_{\gamma_1\gamma_2}A/F_{\gamma_1\gamma_2}^*A.$ }
{\parindent=0pt\v5
{\bf $\Gamma$-filtered module}\par
Let $A$ be a $\Gamma$-filtered $K$-algebra with $\Gamma$-filtration $FA=\{
F_{\gamma}A\}_{\gamma\in\Gamma}$ and $M$ an $A$-module. We say that $M$ is a
$\Gamma$-{\it filtered} $A$-{\it module} if there is a family $FM=\{
F_{\gamma}M\}_{\gamma\in\Gamma}$ consisting of $K$-subspaces $F_{\gamma}M$ of $M$,
such that {\parindent=1.6truecm\vskip 6pt \re{\bf (FM1)}
$M=\cup_{\gamma\in\Gamma}F_{\gamma}M$,\par \re{\bf (FM2)} $F_{\gamma_1}M\subseteq
F_{\gamma_2}M$ if $\gamma_1\preceq\gamma_2$,\par \re{\bf (FM3)}
$F_{\gamma_1}AF_{\gamma_2}M\subseteq F_{\gamma_1\gamma_2}M$,
$\gamma_1,\gamma_2\in\Gamma$.}} \v5 In the definition given above, the family
$FM=\{ F_{\gamma}M\}_{\gamma\in\Gamma}$ is usually called a $\Gamma$-{\it
filtration} of $M$.\par {\parindent=0pt \v5
{\bf Example} (2) Given a $\Gamma$-filtered $K$-algebra $A$ with
$\Gamma$-filtration $FA$, if $\Gamma$ has a smallest element $\gamma_0$ (for
instance $\Gamma =\NZ$), then by the convention we made for $FA$, $1\in
F_{\gamma_0}A$. In this case, any $A$-module $M$ has a $\Gamma$-filtration $FM$. To
see this, let $\{\xi_i\}_{i\in J}$ be a generating set of $M$, i.e., $M=\sum_{i\in
J}A\xi_i$. Put $V=\sum_{i\in J}F_{\gamma_0}A\xi_i$. Then it may be verified
directly that the family $FM=\{ F_{\gamma}M\}_{\gamma\in\Gamma}$ with
$F_{\gamma}M=F_{\gamma}AV$ forms a $\Gamma$-filtration of $M$. \v5
{\bf The associated $\Gamma$-graded module}\par
Let $A$ be a $\Gamma$-filtered $K$-algebra with  $\Gamma$-filtration $FA=\{
F_{\gamma}A\}_{\gamma\in\Gamma}$ and $G(A)$ the associated $\Gamma$-graded
$K$-algebra of $A$. For a $\Gamma$-filtered $A$-module $M$ with
$\Gamma$-filtration $FM=\{ F_{\gamma}M\}_{\gamma\in\Gamma}$, put
$$\begin{array}{l} F_{\gamma}^*M=\bigcup_{\gamma'\prec\gamma}F_{\gamma'}M,\quad \gamma\in\Gamma ,\\
\hbox{where}~F_{\gamma_0}^*M=\{ 0\}~\hbox{if}~A~\hbox{has a smallest element}~\gamma_0.\end{array}$$
The {\it associated $\Gamma$-graded module} of $M$, denoted $G(M)$, is the
$\Gamma$-graded $G(A)$-module defined as
$$G(M)=\bigoplus_{\gamma\in\Gamma}G(M)_{\gamma}~\hbox{with}~G(M)_{\gamma}=F_{\gamma}M/F_{\gamma}^*M,$$
where the module action is given by extending the maps
$$\begin{array}{ccc} G(A)_{\gamma_1}\times G(M)_{\gamma_2}&\mapright{}{}&G(M)_{\gamma_1\gamma_2}\\
(\OV{a_{\gamma_1}},\OV{m_{\gamma_2}})&\mapsto&\OV{a_{\gamma_1}m_{\gamma_2}}\end{array}$$
to $G(A)\times G(M)\r G(M)$, in which $\OV{a_{\gamma_1}}$, $\OV{m_{\gamma_2}}$ are
the images of $a_{\gamma_1}\in F_{\gamma_1}A$, $m_{\gamma_2}\in F_{\gamma_2}A$ in
$G(A)_{\gamma_1}=F_{\gamma_1}A/F_{\gamma_1}^*A$ and
$G(M)_{\gamma_2}=F_{\gamma_2}M/F_{\gamma_2}^*M$ respectively, and
$\OV{a_{\gamma_1}m_{\gamma_2}}$ is the image of $a_{\gamma_1}m_{\gamma_2}\in
F_{\gamma_1\gamma_2}M$ in
$G(M)_{\gamma_1\gamma_2}=F_{\gamma_1\gamma_2}M/F_{\gamma_1\gamma_2}^*M.$ \v5
{\bf $\Gamma$-filtered submodule and induced $\Gamma$-filtration}\par
Let $A$ be a $\Gamma$-filtered $K$-algebra with $\Gamma$-filtration $FA=\{
F_{\gamma}A\}_{\gamma\in\Gamma}$ and $M$ a $\Gamma$-filtered $A$-module with
$\Gamma$-filtration $FM=\{ F_{\gamma}M\}_{\gamma\in\Gamma}$. If $N$ is an
$A$-submodule of $M$ and $N$ has a $\Gamma$-filtration $FN=\{
F_{\gamma}N\}_{\gamma\in\Gamma}$ satisfying $F_{\gamma}N\subseteq F_{\gamma}M$ for
all $\gamma\in\Gamma$, then we call $N$ a $\Gamma$-{\it filtered $A$-submodule} of
$M$. Indeed, any $A$-submodule $N$ of $M$ can be made into a $\Gamma$-filtered
$A$-submodule by using the  {\it induced $\Gamma$-filtration} $FN$ consisting of
$$F_{\gamma}N=N\cap F_{\gamma}M,\quad \gamma\in\Gamma .$$}\par
Furthermore, for an $A$-submodule $N$ of $M$, the quotient $A$-module $M/N$ has the
{\it induced} $\Gamma$-filtration $F(M/N)$ consisting of
$$F_{\gamma}(M/N)=(F_{\gamma}M+N)/N,\quad \gamma\in\Gamma .$$
{\parindent=0pt\par
{\bf $\Gamma$-filtered homomorphism}\par
Let $A$, $B$ be $\Gamma$-filtered $K$-algebras with $FA=\{
F_{\gamma}A\}_{\gamma\in\Gamma}$, $FB=\{ F_{\gamma}B\}_{\gamma\in\Gamma}$,
respectively. A $K$-algebra homomorphism $\varphi$: $A\r B$ is called a
$\Gamma$-{\it filtered $K$-algebra homomorphism} if $\varphi (F_{\gamma}A)\subseteq
F_{\gamma}B$ for all $\gamma\in\Gamma$.}\par Given a $\Gamma$-filtered $K$-algebra
$A$ and two $\Gamma$-filtered $A$-modules $M$, $N$ with $FM=\{
F_{\gamma}M\}_{\gamma\in\Gamma}$, $FN=\{ F_{\gamma}N\}_{\gamma\in\Gamma}$,
respectively. A $\Gamma$-{\it filtered $A$-homomorphism} from $M$ to $N$ is an
$A$-module homomorphism $\psi$: $M\r N$ such that $\psi (F_{\gamma}M)\subseteq
F_{\gamma}N$ for all $\gamma\in\Gamma$.\par A $\Gamma$-filtered $A$-homomorphism
$\psi$: $M\r N$ is said to be {\it strict} if it satisfies
$$\psi (F_{\gamma}M)=\psi (M)\cap F_{\gamma}N,\quad\gamma\in\Gamma .$$
Let $M$ be a $\Gamma$-filtered $A$-module with  $\Gamma$-filtration $FM$ and $N$ a
submodule of $M$. Then, with respect to the filtration $FN=\{ F_{\gamma}N=N\cap
F_{\gamma}M\}_{\gamma\in\Gamma}$ of $N$ and the filtration $F(M/N)=\{
(F_{\gamma}M+N)/N\}_{\gamma\in\Gamma}$ of the quotient module $M/N$, induced by
$FM$, the inclusion map $N\hookrightarrow M$ and the canonical map $M\r M/N$ are
obviously strict $\Gamma$-filtered $A$-homomorphisms. \v5
Verification of the following proposition is an easy exercise. {\parindent=0pt\v5
{\bf 1.1. Proposition} If $\psi$: $M\r N$ is a $\Gamma$-filtered $A$-homomorphism,
then $V=$ Im$\psi$ is a $\Gamma$-filtered $A$-submodule of $N$ with the
$\Gamma$-filtration $FV=\{ F_{\gamma}V=\psi (F_{\gamma}M)\}_{\gamma\in\Gamma}$, and
$W=$ Ker$\psi$ is a $\Gamma$-filtered $A$-submodule of $M$ with the induced
filtration $FW=\{ F_{\gamma}W=W\cap F_{\gamma}M\}_{\gamma\in\Gamma}$.\par\QED \v5
{\bf The associated $\Gamma$-graded $G(A)$-homomorphism}\par
If $\varphi$: $M\r N$ is a $\Gamma$-filtered $A$-homomorphism, then $\varphi$
induces naturally a $\Gamma$-graded $G(A)$-homomorphism:
$$\begin{array}{cccc}
G(\varphi ):&G(M)=\displaystyle{\bigoplus_{\gamma\in\Gamma}}G(M)_{\gamma}&
\mapright{}{}&\displaystyle{\bigoplus_{\gamma\in\Gamma}}G(N)_{\gamma} =G(N)\\
\\
&\sum \OV m&\mapsto&\sum\overline{\varphi (m)}\end{array}$$
where $\OV m$, respectively $\overline{\varphi (m)}$, is the image of $m\in
F_{\gamma}M$ in $G(M)_{\gamma}$, respectively the image of $\varphi (m)$ in
$G(N)_{\gamma}$. \v5 {\bf Remark} If we replace the field $K$ by $\GZ$, then the
text of this section becomes that for $\Gamma$-filtered rings and $\Gamma$-filtered
modules without any modification.} \v5

\section*{2. With $\Gamma$-grading Filtration: $G(R/I)\cong R/\langle\HT (I)\rangle$}
Let $R=\oplus_{\gamma\in\Gamma}R_{\gamma}$ be a $\Gamma$-graded $K$-algebra, where
$\Gamma$ is a totally ordered semigroup with the total ordering $\prec$. In this
section we establish, for an arbitrary ideal $I$ and the quotient algebra $A=R/I$
defined by $I$, the $\Gamma$-graded $K$-algebra isomorphism $G(A)\cong R/\langle\HT
(I)\rangle$ with respect to the $\Gamma$-filtration $FA$ induced by the
$\Gamma$-grading filtration $FR$ as defined in section 1, where $\langle\HT
(I)\rangle$ is the $\Gamma$-graded ideal of $R$ generated by the set of head terms
of  $I$ (see the definition of a head term below). Besides, we conclude a similar
result for an arbitrary left ideal $L$ and the module $M=R/L$. \v5
To start with, note that each element $f\in R$ can be written uniquely as a sum of
finitely many homogeneous elements, say $f=\sum_{i=1}^sr_{\gamma_i}$,
$r_{\gamma_i}\in R_{\gamma_i}$. Assuming
$\gamma_1\succ\gamma_2\succ\cdots\succ\gamma_s$, we define the {\it head term} of
$f$, denoted $\HT (f)$, to be $r_{\gamma_1}$, that is,
$$\HT (f)=r_{\gamma_1} ,$$
and say that $f$ is of {\it degree} $\gamma_1$, denoted $d(f)=\gamma_1$. For a
subset $S\subset R$, we define the set of head terms of $S$ as
$$\HT (S)=\left\{ \HT (r)~\Big |~r\in S\right\} .$$\par
Let $I$ be an ideal of $R$. As $\HT (I)$ consists of homogeneous elements, the
ideal $\langle\HT (I)\rangle$ of $R$ is $\Gamma$-graded, and hence, the quotient
algebra $\OV A=R/\langle\HT (I)\rangle$ is a $\Gamma$-graded $K$-algebra with the
$\Gamma$-gradation $\{ \OV A_{\gamma}=(R_{\gamma}+\langle\HT (I)\rangle
)/\langle\HT (I)\rangle\}_{\gamma\in\Gamma}$. Consider the $\Gamma$-grading
filtration $FR=\{ F_{\gamma}R\}_{\gamma\in\Gamma}$ of $R$ in the sense of section
1, where
$$F_{\gamma}R=\bigoplus_{\gamma '\preceq\gamma}F_{\gamma '}R,\quad \gamma\in\Gamma .$$
Then the  quotient algebra $A=R/I$ has the induced $\Gamma$-filtration $FA=\{
F_{\gamma}A=(F_{\gamma}R+I)/I\}_{\gamma\in\Gamma}$. Taking the associated
$\Gamma$-graded $K$-algebra $G(A)$ of $A$ defined by  $FA$ into account, we have
the following key result of this paper. {\parindent=0pt\v5
{\bf 2.1. Theorem} With notation as fixed above, we have a $\Gamma$-graded
$K$-algebra isomorphism
$$G(A)\cong \OV A=R/\langle\HT (I)\rangle .$$
{\bf Proof} First, recall that $\Gamma$ is ordered by the total ordering $\prec$. By
the definition of $G(A)$, for $\gamma\in\Gamma$,
$G(A)_{\gamma}=F_{\gamma}A/F_{\gamma}^*A$ with $F_{\gamma}A=(F_{\gamma}R+I)/I$ and,
as a $K$-subspace,
$$\begin{array}{rcl} F_{\gamma}^*A=\displaystyle{\bigcup_{\gamma '\prec\gamma}}F_{\gamma '}A&=&
\displaystyle{\bigcup_{\gamma '\prec\gamma}}\FRAC{F_{\gamma '}R+I}{I}\\
\\
&=&\FRAC{\bigcup_{\gamma '\prec\gamma}F_{\gamma
'}R+I}{I}=\FRAC{F_{\gamma}^*R+I}{I}.\end{array}$$ It turns out that there are
canonical isomorphisms of $K$-subspaces
$$\FRAC{R_{\gamma}\oplus F_{\gamma}^*R}{(I\cap F_{\gamma}R)+F_{\gamma}^*R}=
\FRAC{F_{\gamma}R}{(I\cap
F_{\gamma}R)+F_{\gamma}^*R}~\mapright{\cong}{}~G(A)_{\gamma},\quad \gamma\in\Gamma ,
$$ and consequently, we can extend the natural epimorphisms of $K$-subspaces
$$\phi_{\gamma} :~ R_{\gamma}~\mapright{}{}~\FRAC{R_{\gamma}\oplus F_{\gamma}^*R}{(I\cap F_{\gamma}R)+F_{\gamma}^*R},\quad\gamma\in\Gamma ,$$
to define a $\Gamma$-graded $K$-algebra epimorphism
$$\phi :~R~\mapright{}{}~G(A).$$
We claim that Ker$\phi =\langle\HT (I)\rangle$. To see this, noticing $\langle\HT
(I)\rangle$ is a $\Gamma$-graded ideal, it is sufficient to prove the equalities
$$\hbox{Ker}\phi_{\gamma} =\langle\HT (I)\rangle\cap R_{\gamma},\quad \gamma\in\Gamma .$$
Suppose $r_{\gamma}\in$ Ker$\phi_{\gamma}\subset R_{\gamma}$. Then $r_{\gamma}\in$
$(I\cap F_{\gamma}R)+F_{\gamma}^*R$. If $r_{\gamma}\ne 0$, then, as a homogeneous
element of degree $\gamma$, $r_{\gamma}=\HT (f)$ for some $f\in I\cap F_{\gamma}R$.
This shows that $r_{\gamma}\in\langle\HT (I)\rangle\cap R_{\gamma}$. Hence
Ker$\phi_{\gamma}\subseteq\langle\HT (I)\rangle\cap R_{\gamma}$. Conversely,
suppose $r_{\gamma}\in$ $\langle\HT (I)\rangle\cap R_{\gamma}$. Then, as a
homogeneous element of degree $\gamma$, $r_{\gamma}=\sum_{i=1}^s v_i\HT (f_i)w_i$,
where $v_i$, $w_i$ are homogeneous elements of $R$ and $f_i\in I$. Write $f_i=\HT
(f_i)+f_i'$ such that $d(f_i')\prec d(f_i)$, $i=1,...,s$.   By the fact that
$\Gamma$ is an ordered semigroup with the total ordering $\prec$, we may see that
the expression
$$r_{\gamma}=\sum^s_{i=1}v_if_iw_i-\sum^s_{i=1}v_if_i'w_i$$
satisfies $\sum^s_{i=1}v_if_iw_i\in I\cap F_{\gamma}R$ and
$\sum^s_{i=1}v_if_i'w_i\in F_{\gamma}^*R$. This shows that $r_{\gamma}\in (I\cap
F_{\gamma}R)+F_{\gamma}^*R$, i.e., $r_{\gamma}\in$ Ker$\phi_{\gamma}$. Hence,
$\langle\HT (I)\rangle\cap R_{\gamma}\subseteq$ Ker$\phi_{\gamma}$. Summing up, we
conclude the desired equalities Ker$\phi_{\gamma}=\langle\HT (I)\rangle\cap
R_{\gamma}$, $\gamma\in\Gamma$.\QED \v5
{\bf Remark} Obviously, if $I$ is a $\Gamma$-graded ideal of $R$, then $A=R/I=G(A)$
with respect to $FA$ induced by the $\Gamma$-grading filtration $FR$ of $R$. }\v5
We illustrate Theorem 2.1 by two classical examples.{\parindent=0pt\v5
{\bf Example} (1) Let $\textsf{g}$ be a $K$-Lie algebra with the $K$-basis $\{
x_1,...,x_n\}$ and the bracket product
$$[x_i,~x_j]=\sum_{\ell =1}^n\lambda_{ij}^{\ell}x_{\ell},\quad 1\le i<j\le n,
~\lambda_{ij}^{\ell}\in K ,$$  and let $U(\textsf{g})$ be the universal enveloping
algebra of $\textsf{g}$. Then  $U(\textsf{g})=\KS /I$, where $\KS =K\langle
X_1,...,X_n\rangle$ is the free $K$-algebra generated by $X=\{ X_1,...,X_n\}$ over
$K$, and the ideal $I$ is generated by $\G =\{ X_jX_i-X_iX_j-\sum^n_{\ell
=1}\lambda_{ij}^nX_{\ell}~|~ 1\le i<j\le n\}$. If we consider the
$\mathbb{N}$-grading filtration $F\KS$ of $\KS$ as defined in section 1, then
$F\KS$ induces the natural $\mathbb{N}$-filtration $FU(\textsf{g})=\{ (F_p\KS
+I)/I\}_{p\in\mathbb{N}}$ of $U(\textsf{g})$. Hence, by Theorem 2.1,
$U(\textsf{g})$ has the associated $\mathbb{N}$-graded algebra
$G(U(\textsf{g}))\cong\KS/\langle\HT (I)\rangle$. It is well-known that $\G$ is a
Gr\"obner basis for the ideal $I$ with respect to a graded monomial ordering on
$\KS$ (see [Mor]). By [Li1], $\HT (\G)$ is a Gr\"obner basis of $\langle \HT
(I)\rangle$ (also see later section 5). Hence, $G(U(\textsf{g}))\cong
K[x_1,...,x_n]$, the commutative polynomial $K$-algebra in $n$ variables. So, the
classical PBW theorem is recaptured.  \vskip 6pt
(2) Let $A_n(K)$ be the $n$-th Weyl algebra, that is, $A_n(K)\cong K\langle
X,~Y\rangle /I$, where $K\langle X,~Y\rangle=K\langle
X_1,...,X_n,Y_1,...,Y_n\rangle$ is the free $K$-algebra generated by $X\cup Y=\{
X_1,...,X_n,Y_1,...,Y_n\}$ over $K$, and the ideal $I$ is generated by the set $\G$
consisting of
$$\begin{array}{ll} X_iX_j-X_jX_i,~Y_iY_j-Y_jY_i,&1\le i<j\le n,\\
Y_jX_i-X_iY_j-\delta_{ij}~\hbox{the Kronecker delta},&1\le i,j\le n.\end{array}$$
If we consider the $\mathbb{N}$-grading filtration $FK\langle X,~Y\rangle$ of
$K\langle X,~Y\rangle$ as defined in section 1, then $FK\langle X,~Y\rangle$
induces the natural $\mathbb{N}$-filtration $F_pA_n(K)=\{ (F_p\KS
+I)/I\}_{p\in\mathbb{N}}$ of $A_n(K)$. Hence, by Theorem 2.1, $A_n(K)$ has the
associated $\mathbb{N}$-graded algebra $G(A_n(K))\cong K\langle
X,~Y\rangle/\langle\HT (I)\rangle$. It is equally well-known that $\G$ is a
Gr\"obner basis for the ideal $I$ with respect to a graded monomial ordering on
$K\langle X,~Y\rangle$ (see [Mor]). By [Li1], $\HT (\G)$ is a Gr\"obner basis of
$\langle \HT (I)\rangle$ (also see later section 5). Hence, the classical result
$G(A_n(K))\cong K[x_1,...,x_n,y_1,...,y_n]$ is also recaptured, where the latter is
the commutative polynomial $K$-algebra in $2n$ variables. } \v5
In general, suppose the ideal $I$ is generated by $F=\{ f_i\}_{i\in J}$. Put $\HT
(F)=\{ \HT (f_i)~|~f_i\in F\}$. Then, naturally, we expect that the equality
$$\langle\HT (I)\rangle =\langle\HT (F)\rangle \leqno{(*)}$$
holds, and consequently we would have $G(A)\cong R/\langle\HT (F)\rangle$. In other
words, the equality $(*)$ amounts to propose a general version of PBW theorem. To
realize  this property effectively in later section 5, let us examine how a
generating set of $I$ gives rise to a generating set for $\langle\HT (I)\rangle$,
and vice versa.{\parindent=0pt\v5
{\bf 2.2. Proposition} Let $F=\{ f_i\}_{i\in J}$ be a subset of the ideal $I$. The
following two statements hold.\par
(i) Suppose that $F$ is a generating set of $I$ having the property that each
nonzero $f\in I$ has a presentation
$$\begin{array}{rcl} f&=&\sum v_jf_jw_j,~\hbox{in which}~v_j,~w_j~\hbox{are homogeneous elements of}~R,~f_j\in F,\\
&~~~&\hbox{such that}~v_jf_jw_j\ne 0~\hbox{and}~d(\HT (v_j\HT (f_j)w_j))\preceq
d(f),\\
&~~~&\hbox{where}~f_j~\hbox{may appear repeatedly}.\end{array}$$
Then $\langle\HT (I)\rangle =\langle\HT (F)\rangle$.\par
(ii) In the case that $\prec$ is a well-ordering on $\Gamma$, if $\langle\HT
(I)\rangle =\langle\HT (F)\rangle$, then $F$ is a generating set of $I$ having the
property mentioned in (i) above.\vskip 6pt
{\bf Proof} (i)  By definition, if $f\in R$, $f\ne 0$ and $\HT (f)\in R_{\gamma}$,
then $d(f)=d(\HT (f))=\gamma$. Since $\Gamma$ is an ordered semigroup with the
total ordering $\prec$, by the assumption on the presentation $f=\sum v_jf_jw_j$,
the head term $\HT (f)$ of $f$ must have the form $\HT (f)=\sum v_j\HT (f_j')w_j$,
i.e., $\HT (f)\in\langle\HT (F)\rangle$. Hence $\langle\HT (I)\rangle =\langle\HT
(F)\rangle$.\par
(ii)  For $f\in I$, $f\ne 0$, suppose $\HT (f)\in R_{\gamma}$. By the assumption,
the homogeneous element $\HT (f)$ can be written as
$$\HT (f)=\sum^s_{j=1}v_j\HT (f_j)w_j,$$
in which $v_j$, $w_j$ are homogeneous of $R$ and  $f_j\in F$, satisfying $v_j\HT
(f_j)w_j\ne 0$ and $d(v_j\HT (f_j)w_j)=d(\HT (f))=d(f)=\gamma$, $j=1,...,s$, where
$\HT (f_j)$ may appear repeatedly.  Thus, writing each $f_j$ as $f_j=\HT
(f_j)+f_j'$ such that $d(f_j')\prec d(f)$,  we have
$$\HT (f)=\sum^s_{j=1}v_jf_jw_j-\sum^s_{j=1}v_jf_j'w_j, $$
in which each $v_jf_jw_j\ne 0$, $d(\sum^s_{j=1}v_jf_jw_j)=d(\HT (f))=d(f)=\gamma$,
and $d(\sum^s_{j=1}v_jf_j'w_j)\prec\gamma$.  It turns out that the element
$$f'=f-\sum^s_{j=1}v_jf_jw_j\in I$$
has $d(f')\prec d(f)=\gamma$. For $f'$, we may repeat the same procedure and get
some
$$f''=f'-\sum^m_{k=1}v_kf_kw_k\in I$$
with $d(f'')\prec\gamma '$, where $v_k$, $w_k$ are homogeneous elements of $R$,
$f_k\in F$,  satisfying each $v_kf_kw_k\ne 0$ and $d(\sum^m_{k=1}v_kf_kw_k)=\gamma
'$. Since $\prec$ is a well-ordering, after a finite number of repetitions, such
reduction procedure of decreasing degrees must stop to give us an expression
$f=\sum v_jf_jw_j$ with the desired property.\QED} \v5
We finish this section by remarking that the proof of Theorem 2.1 and Proposition
2.2 may be carried to deal with a left ideal $L$ of $R$ and the module $M=R/L$
directly so long as we replace $I$ by $L$ and consider only left-hand side action.
We mention the result below but will not dig in detail on this topic in this
paper.{\parindent=0pt\v5

{\bf 2.3. Theorem} Let $L$ be an arbitrary left ideal of $R$ and $M=R/L$ the
quotient module of $R$ determined by $L$. Considering the $\Gamma$-filtration
$FL=\{ F_{\gamma}L=L\cap F_{\gamma}R\}_{\gamma\in\Gamma}$ of $L$ and the
$\Gamma$-filtration $FM=\{ F_{\gamma}M=(F_{\gamma}R+L)/L\}_{\gamma\in\Gamma}$ of
$M$, induced by the $\Gamma$-grading filtration $FR$ of $R$, let $G(M)$ be the
associated $\Gamma$-graded $G(A)$-module. Then we have a $\Gamma$-graded
$R$-isomorphism
$$G(M)\cong R/\langle\HT (L)],$$
where $\langle\HT (L)]$ denotes the $\Gamma$-graded left ideal of $R$ generated by
$\HT (L)$.\par\QED \v5
{\bf 2.4. Proposition}  Let $F=\{ f_i\}_{i\in J}$ be a subset of a left ideal $L$
of $R$. The following two statements hold.\par
(i) Suppose that $F$ is a generating set of $L$ having the property that each
nonzero $f\in L$ has a presentation
$$\begin{array}{rcl} f&=&\sum v_jf_j,~\hbox{in which the}~v_j\hbox{are homogeneous elements of}~R,~f_j\in F,\\
&~~~&\hbox{such that}~v_jf_j\ne 0~\hbox{and}~d(\HT (v_j\HT (f_j)))\preceq d(f),\\
&~~~&\hbox{where}~f_j~\hbox{may appear repeatedly}.\end{array}$$
Then $\langle\HT (I)] =\langle\HT (F)]$.\par
(ii) In the case that $\prec$ is a well-ordering on $\Gamma$, if $\langle\HT (I)]
=\langle\HT (F)]$, then $F$ is a generating set of $L$ having the property
mentioned in (i) above.\par\QED}\v5

\section*{3. Basic Lifting Properties}
In this and the next section we explore the structural relation between
$\Gamma$-filtered $A$-modules and $\Gamma$-graded $G(A)$-modules, where $\Gamma$ is
an ordered monoid with a well-ordering, that leads to many nice lifting properties.
In view of section 2, these lifting properties  provide us with a firm basis to
study quotient algebras $A=R/I$ of a $\Gamma$-graded $K$-algebra $R$ via  the
quotient algebra $R/\langle\HT (I)\rangle=G(A)$ with respect to the
$\Gamma$-filtration $FA$ induced by the $\Gamma$-grading filtration $FR$ of $R$,
especially when $R$ is taken to be a free algebra, or a commutative polynomial
algebra, or a path algebra, or some other commonly used graded algebra such as the
coordinate rings of quantum affine $K$-spaces, Gr\"obner bases may be used to
realize the lifting properties effectively (see later sections 5 -- 7). \v5
Let $\Gamma$ be an ordered {\it monoid} with the {\it well-ordering} $\prec$. If
$\gamma_0$ is the identity element of $\Gamma$, we assume that $\gamma_0$ is the
{\it smallest} element in $\Gamma$. All notations used in previous sections are
maintained. \par
Let $A$ be a $\Gamma$-filtered $K$-algebra with  $\Gamma$-filtration $FA$, and let
$G(A)$ be the associated $\Gamma$-graded $K$-algebra of $A$. Then $1\in
F_{\gamma_0}A$ by the assumption on $\Gamma$ made above and the convention fixed in
the definition of $FA$ (section 1). \v5
Let $M$ be a $\Gamma$-filtered $A$-module with  $\Gamma$-filtration $FM=\{
F_{\gamma}M\}_{\gamma\in\Gamma}$, and let $G(M)$ be the associated $\Gamma$-graded
$G(A)$-module of $M$, in the sense of section 3. Since $\prec$ is a well-ordering
on $\Gamma$, for each element $m\in M$, we define the {\it degree} of $m$, denoted
$d(m)$, as
$$d(m)=\min\left\{ \gamma\in\Gamma ~\Big |~m\in F_{\gamma}M\right\} .$$
If $m\ne 0$ and $d(m)=\gamma$, then we write $\sigma (m)$ for the corresponding {\it
nonzero} homogeneous element of degree $\gamma$ in
$G(M)_{\gamma}=F_{\gamma}M/F_{\gamma}^*M$. \v5 As to the $\sigma$-element defined
above, We first note an easily verified but useful property:
$$\forall~a\in A,~m\in M,~\hbox{either}~\sigma (a)\sigma (m)=0~\hbox{or}~\sigma (a)\sigma (m)=\sigma (am).\leqno{(\sigma )}$$\par
We deal first with $K$-basis, divisors of zeros and primeness (semi-primeness).
Recall that a ring $S$ is a domain if $S$ does not have divisors of zero; and $S$
is called a prime (semi-prime) ring if $s_1Ss_2\ne 0$ for any nonzero $s_1,s_2\in
S$ (if $sSs\ne 0$ for any nonzero $s\in S$). {\parindent=0pt\v5
{\bf 3.1. Theorem} Let $A$ be an arbitrary $\Gamma$-filtered $K$-algebra with
$\Gamma$-filtration $FA$, and let $G(A)$ be the associated $\Gamma$-graded
$K$-algebra of $A$. The following statements hold, especially when $A=R/I$ and
$G(A)=R/\langle\HT (I)\rangle$ as in Theorem 2.1.\par
(i) Suppose that $\{ a_i\}_{i\in J}$ is a subset of $A$ such that $\{ \sigma
(a_i)\}_{i\in J}$ forms a $K$-basis for $G(A)$, then $\{ a_i\}_{i\in J}$ is a
$K$-basis of $A$. Hence, if $G(A)$ is finite dimensional over $K$ then so is
$A$.\par
(ii) If $G(A)$ is a domain then so is $A$.\par
(iii) If $G(A)$ is a prime (semi-prime) ring then so is $A$. \vskip 6pt {\bf Proof}
(i) We show first that $\{ a_i\}_{i\in J}$ is $K$-linearly independent, namely,  if
$a=\sum_{j=1}^s\lambda_{i_j}a_{i_j}=0$, where $\lambda_i\in K$ and $a_{i_j}\in\{
a_i\}_{i\in J}$, then all coefficients $\lambda_{i_j}=0$. To this end, assume that
$a_{i_1}$, $a_{i_2}$, ..., $a_{i_t}$, $t\le s$, all have the same degree $\gamma$
that is the highest degree among all terms in the linear expression of $a$. Then
since $K\subseteq F_{\gamma_0}A$, taking the image of $a$ in
$G(A)_{\gamma}=F_{\gamma}A/F_{\gamma}^*A$ into account, we have
$$\lambda_{i_1}\sigma (a_{i_1})+\lambda_{i_2}\sigma (a_{i_2})+
\cdots \lambda_{i_t}\sigma (a_{i_t})=0$$
in $G(A)_{\gamma}$. By the $K$-linear independence of $\{\sigma (a_i)\}_{i\in J}$, $\lambda_{i_1}=
\lambda_{i_2}=\cdots =\lambda_{i_t}=0$. Similarly we assert that all other coefficients
in the linear expression of $a$ are equal to 0. This proves the $K$-linear independence of $\{a_i\}_{i\in J}$. }\par
Next, we conclude that $A$, as a $K$-space, is spanned by $\{ a_i\}_{i\in J}$. To see
this, let $a\in F_{\gamma}A-F_{\gamma}^*A$, i.e., $d(a)=\gamma$. Then, by the assumption,
$\sigma (a)$ can be written uniquely as a linear combination of $\sigma (a_i)$'s, say
$$\sigma (a)=\sum_{j=1}^s\lambda_{i_j}\sigma (a_{i_j}),\quad\lambda_{i_j}\in K,~a_{i_j}\in\{ a_i\}_{i\in J} .$$
As $\sigma (a)$ is a homogeneous element of degree $\gamma$ in $G(A)$ and $K\subseteq G(A)_{\gamma_0}=F_{\gamma_0}A$, it follows that
all homogeneous elements $\sigma (a_{i_j})$ in the linear expression of $\sigma (a)$ have the same degree
$\gamma$. Thus, by the definition of a $\sigma$-element, we have
$$a'=a-\sum^s_{j=1}\lambda_{i_j}a_{i_j}\in F_{\gamma}^*A.$$
Suppose $a'\in F_{\gamma '}A-F_{\gamma '}^*A$, i.e., $d(a')=\gamma '\prec\gamma$.
By a similar argument we may get $a''=a'-\sum_{\ell
=1}^m\lambda_{i_{\ell}}a_{i_{\ell}}\in F_{\gamma '}^*A$ with $\lambda_{i_{\ell}}\in
K$ and $a_{i_{\ell}}\in\{ a_i\}_{i\in J}$. If $d(a'')=\gamma ''$, then
$\gamma\succ \gamma '\succ\gamma ''$. Since $\prec$ is a well-ordering, after
repeating the procedure of reducing degrees for a finite number of steps we must
have $a\in K$-span$\{ a_i\}_{i\in J}$, as desired. {\parindent=0pt\par (ii) Let
$a,b\in A$ be
 nonzero elements of degree $\gamma_1$, $\gamma_2$ respectively. Then $\sigma (a)$,
 $\sigma (b)$ are nonzero homogeneous elements of $G(A)$ and so $\sigma (a)\sigma
 (b)= \sigma (ab)\ne 0$. This means $ab\not\in F_{\gamma_1\gamma_2}^*A$, and hence
 $ab\ne 0$.\par
(iii) If $a,b\in A$ are nonzero, then $\sigma (a)$, $\sigma (b)$ are
 nonzero homogeneous elements of $G(A)$ and so $\sigma (a)G(A)\sigma (b)\ne \{ 0\}$.
 It follows that there is a homogeneous element $\sigma (c)\in G(A)$, represented by
 $c\in A$, such that $\sigma (a)\sigma (c)\sigma (b)\ne 0$. But this means $acb\ne
 0$. Hence $aAb\ne \{ 0\}$, i.e., $A$ is prime. A similar argument holds for the
 semi-primeness.\QED} \v5
Next, we focus on modules. {\parindent=0pt\v5
{\bf  3.2. Proposition} Let $M$ be an $A$-module.\par (i) Let $FM=\{
 F_{\gamma}M\}_{\gamma\in\Gamma}$ be a $\Gamma$-filtration of $M$. If
 $G(M)=\sum_{i\in J}G(A)\sigma (\xi_i)$ with $\xi_i\in M$ and $d
 (\xi_i)=\gamma_i\in\Gamma$, then $M=\sum_{i\in J}A\xi_i$ with
$$F_{\gamma}M=\sum_{i\in J}\left (\sum_{s\preceq\gamma ,~s_i\gamma_i=s}F_{s_i}A\right )\xi_i,\quad \gamma\in\Gamma .$$
In particular, if $G(M)$ is finitely generated then so is $M$.\par
(ii) If $M$ is finitely generated, then $M$ has a $\Gamma$-filtration $FM=\{
F_{\gamma}M\}_{\gamma\in\Gamma}$ such that $G(M)$ is a finitely generated
$G(A)$-module. Indeed, if $M=\sum_{i=1}^nA\xi_i$ and $\{ \xi_1,...,\xi_n\}$ is a
{\it minimal} set of generators for $M$, then the desired $\Gamma$-filtration $FM$
consists of
$$F_{\gamma}M=\displaystyle{\sum^n_{i=1}}\left (
\displaystyle{\sum_{s\preceq\gamma ,~s_i\gamma_i=s}}F_{s_i}A \right ) \xi_i,\quad
\gamma\in\Gamma,$$ where $\gamma_1,...,\gamma_n\in\Gamma$ are chosen arbitrarily.
\vskip 6pt                                                         {\bf Proof} (i)
Since $G(M)=\sum_{i\in J}G(A)\sigma (\xi_i)$ with $\xi_i\in M$ and
$d(\xi_i)=\gamma_i\in\Gamma$, we have
$$G(M)_{\gamma}=\displaystyle{\sum_{i\in J,~\rho_i\gamma_i=\gamma}}G(A)_{\rho_i}\sigma (\xi_i),\quad \gamma\in\Gamma .$$
Hence, for any $m\in F_{\gamma}M$, $m=\sum a_{\rho_i}\xi_i+m'$, where $a_{\rho_i}\in
F_{\rho_i}A$, $\rho_i\gamma_i=\gamma$, and $m'\in F_{\gamma}^*M$. Assume
$d(m')=\gamma '$. Then, similarly we have $m'=\sum a_{\mu_i}\xi_i+m''$, where
$a_{\mu_i}\in F_{\mu_i}A$, $\mu_i\gamma_i=\gamma '$, and $m''\in F_{\gamma '}^*M$.
Suppose $d(m '')=\gamma ''$. Then, $\gamma\succ\gamma '\succ\gamma ''$. As $\prec$
is a well-ordering, after repeating the procedure of reducing degrees
for a finite number of steps,  we should arrive at
$$m\in\displaystyle{\sum_{i\in J}}\left (
\displaystyle{\sum_{s\le\gamma ,~s_i\gamma_i=s}}F_{s_i}A \right ) \xi_i.$$ Since
$m$ is taken arbitrarily, this shows that
$$F_uM=\displaystyle{\sum_{i\in J}}\left (
\displaystyle{\sum_{s\le\gamma ,~s_i\gamma_i=s}}F_{s_i}A \right ) \xi_i,$$ and
therefore $M=\sum_{i\in J}A\xi_i$. \par (ii) Suppose $M=\sum^n_{i=1}A\xi_i$ and $\{
\xi_1,...,\xi_n\}$ is a {\it minimal} set of generators for $M$. Choose
$\gamma_1,...,\gamma_n\in\Gamma$ arbitrarily. Then, since each $m\in M$ has a
presentation $m=\sum_{i=1}^na_i\xi_i$ with $a_i\in F_{s_i}A-F_{s_i}^*A$ for some
$s_i\in\Gamma$, it is easy to see that the $K$-subspaces
$$F_{\gamma}M=\displaystyle{\sum^n_{i=1}}\left (
\displaystyle{\sum_{s\le\gamma ,~s_i\gamma_i=s}}F_{s_i}A \right ) \xi_i,\quad
\gamma\in\Gamma,$$ form a $\Gamma$-filtration $FM=\{
F_{\gamma}M\}_{\gamma\in\Gamma}$ for $M$. Furthermore, note that $1\in
F_{\gamma_0}A$, where $\gamma_0$ is the identity element of $\Gamma$. It follows
from the construction of $FM$ and the minimality of $\{\xi_1,...,\xi_n\}$ (as a set
of generators of $M$) that $\xi_i\in F_{\gamma_i}M-F_{\gamma_i}^*M$, i.e.,
$d(\xi_i)=\gamma_i$, $i=1,...,n$. Thus, by the foregoing property $(\sigma )$ of
the associated $\sigma$-elements, it is not difficult to verify that
$$G(M)_{\gamma}=F_{\gamma}M/F_{\gamma}^*M=\displaystyle{\sum_{1\le i\le n,~s_i\gamma_i=\gamma}}G(A)_{s_i}\sigma (\xi_i), \quad \gamma\in\Gamma ,$$
and hence $G(M)=\oplus_{\gamma\in\Gamma}G(M)_{\gamma}=\sum^n_{i=1}G(A)\sigma
(\xi_i)$.\QED} \v5
Recall that a sequence
$$\cdots~\mapright{\varphi_{i-2}}{}~M_{i-1}~\mapright{\varphi_{i-1}}{}~M_{i}~\mapright{\varphi_i}{}~M_{i+1}~\mapright{\varphi_{i+1}}{}~\cdots$$
of $A$-modules and $A$-homomorphisms is said to be {\it exact} if Ker$\varphi_k =$
Im$\varphi_{k-1}$ holds for every $k$.\par
$\Gamma$-filtered homomorphisms and the associated $\Gamma$-graded homomorphisms
considered below are in the sense of section 1. {\parindent=0pt\v5
{\bf 3.3. Proposition} Let
$$L~\mapright{\varphi}{}~M~\mapright{\psi}{}~N\leqno{(*)}$$
be a sequence of $\Gamma$-filtered $A$-modules and $\Gamma$-filtered
$A$-homomorphisms satisfying $\psi\circ\varphi =0$. Then the following properties
are equivalent.\par
(i) The sequence $(*)$ is exact and $\varphi$ and $\psi$ are strict.\par
(ii) The associated sequence of $\Gamma$-graded $G(A)$-modules and $\Gamma$-graded
$G(A)$-homomorphisms
$$G(L)~\mapright{G(\varphi)}{}~G(M)~\mapright{G(\psi )}{}~G(N)\leqno{G(*)}$$
is exact . \vskip 6pt
{\bf Proof} For an element $x\in F_{\gamma}M$, we use $\OV x$ to denote the image
of $x$ in $F_{\gamma}M/F_{\gamma}^*M=G(M)_{\gamma}$. Similar notation is used for
elements in $L$ and $N$.\par (i) $\Rightarrow$ (ii) By the definition of the
associated $\Gamma$-graded $G(A)$-homomorphism of a $\Gamma$-filtered
$A$-homomorphism, it is clear that Im$G(\varphi )\subseteq$ Ker$G(\psi )$. To prove
the converse inclusion, for $m\in F_{\gamma}M-F_{\gamma}^*M$, i.e., $d(m)=\gamma$,
suppose $0=G(\psi )(\OV m)=\OV{\psi (m)}$. If $\psi (m)=0$, then $m\in$ Ker$\psi =$
Im$\varphi$ and there is some $\ell\in L$ such that
$$m=\varphi (\ell )\in\varphi (L)\cap F_{\gamma}M=\varphi (F_{\gamma}L).$$
Obviously we may assume $\ell\in F_{\gamma}L$, and thus, $\OV m=\OV{\varphi (\ell
)}=G(\varphi )(\OV{\ell})$, i.e., $\OV m\in$ Im$G(\varphi )$. If $\psi (m)\ne 0$,
then since $0=G(\psi )(\OV m)=\OV{\psi (m)}\in G(N)_{\gamma}$, we have $\psi (m)\in
F_{\gamma '}N-F_{\gamma '}^*N$ for some $\gamma '\prec\gamma$, i.e., $\psi (m)\in
\psi (M)\cap F_{\gamma '}N=\psi (F_{\gamma '}M)$. This yields $\psi (m)=\psi (m')$
for some $m'\in F_{\gamma '}M$, and hence
$$m-m'\in~\hbox{Ker}\psi\cap F_{\gamma}M=\varphi (L)\cap F_{\gamma}M=\varphi (F_{\gamma}L).$$
Let $m-m'=\varphi (\ell ')$ for some $\ell '\in F_{\gamma}L$. Then $\OV
m=\OV{m-m'}=\OV{\varphi (\ell ')}=G(\varphi )(\OV{\ell '})$. this shows that $\OV
m\in$ Im$G(\varphi )$. As $m$ is taken arbitrarily, so we have Ker$G(\psi
)\subseteq$ Im$G(\varphi )$. Therefore, we conclude Ker$G(\psi )=$ Im$G(\varphi )$,
that is, the sequence $G(*)$ is exact.\par (ii) $\Rightarrow$ (i)  Suppose that the
graded sequence $G(*)$ is exact. Let us show that the sequence $(*)$ is exact first.
If $\psi (m)=0$ with $m\in F_{\gamma}M-F_{\gamma}^*M$, i.e., $d(m)=\gamma$, then
$G(\psi )(\OV m)=0$ with $\OV m=\sigma (m)\in G(M)_{\gamma}$. It follows that $\OV
m=G(\varphi )(\overline{\ell '})=\overline{\varphi (\ell ')}$ for some $\ell '\in
F_{\gamma}L-F_{\gamma}^*L$. Hence $m-\varphi (\ell ')=m'$ for some $m'\in F_{\gamma
'}M$ with $\gamma '\prec\gamma$. Thus, $\psi (m')=\psi (m-\varphi (\ell '))=0$.
Similarly, $m'-\varphi (\ell '')=m''$ with $\ell ''\in L$ and $m''\in F_{\gamma
''}M$ with $\gamma ''\prec\gamma '$. As the chain
$$\gamma\succ \gamma '\succ\gamma ''\succ\cdots$$
cannot be infinite, for $\prec$ is a well-ordering, after repeating the reduction
procedure for a finite number of steps we arrive at $m=\varphi (\ell)$ for some
$\ell\in L$. This shows that Ker$\psi\subseteq\varphi (L)$. Therefore, Ker$\psi
=\varphi (L)$ and the exactness of the sequence $(*)$ follows.}\par
As to the strictness of $\varphi$ and $\psi$, we prove it only for $\psi$ because a
similar argument is valid for $\varphi$. Let $f\in F_{\gamma}N\cap\psi (M)$ and
$f\not\in F_{\gamma}^*N$. Then $f=\psi (m)$ for some $m\in F_wM$. Suppose
$\gamma\preceq w$. If $w =\gamma$, then $f=\psi (m)\in\psi (F_{\gamma}M)$. If
$\gamma\prec w$, then since $f\in F_{\gamma}N$, we have $G(\psi)(\OV
m)=\overline{\psi (m)}=0$ in $G(N)$. By the exactness, $\OV m=G(\varphi
)(\OV{\ell}) =\overline{\varphi (\ell )}$ for some $\ell\in F_wL$. Put
$m'=m-\varphi (\ell )$. Then $m'\in F_{w'}M$ with $w'\prec w$, and $\psi (m')=\psi
(m-\varphi (\ell ))=\psi (m)=f$. If $\gamma\prec w'$, then similarly we may find
$m''\in F_{w''}M$ with $w''\prec w'$, such that  $\psi (m'')=f$.  Note that the
chain
$$w\succ w'\succ w''\succ\cdots\succ\gamma$$
has finite length in $\Gamma$. So the reduction procedure above stops after a
finite number of steps, and eventually we have  $f=\psi (m_{\gamma})\in\psi
(F_{\gamma}M)$. This shows that $F_{\gamma}N\cap\psi (M)\subset \psi
(F_{\gamma}M)$, that is, $\psi$ is strict.\QED {\parindent=0pt\v5
{\bf 3.4. Corollary} (i) Let $\varphi$: $M\r N$ be a $\Gamma$-filtered
$A$-homomorphism. Then $G(\varphi )$ is injective, respectively surjective, if and
only if $\varphi$ is injective, respectively surjective, and $\varphi$ is
strict.\par (ii) Let $N$, $W$ be submodules of a $\Gamma$-filtered $A$-module $M$
with  $\Gamma$-filtration $FM=\{ F_{\gamma}M\}_{\gamma\in\Gamma}$. Consider the
$\Gamma$-filtration $FN=\{ F_{\gamma}N=N\cap F_{\gamma}M\}_{\gamma\in\Gamma}$ of
$N$ and the $\Gamma$-filtration $FW=\{ F_{\gamma}W=W\cap
F_{\gamma}M\}_{\gamma\in\Gamma}$, induced $FM$ respectively.  If $N\subseteq W$,
then $G(N)\subseteq G(W)$; and if $G(N)=G(W)$ then $N=W$.\par\QED} \v5
We summarize some immediate applications of previous results in the following
theorem. {\parindent=0pt\v5
{\bf 3.5. Theorem} Let $A$ be an arbitrary $\Gamma$-filtered $K$-algebra with
$\Gamma$-filtration $FA$, and let $G(A)$ be the associated $\Gamma$-graded
$K$-algebra of $A$.  The following statements hold, especially when $A=R/I$ and
$G(A)=R/\langle\HT (I)\rangle$ as in Theorem 2.1.\par
(i) Suppose that $G(A)$ is $\Gamma$-graded left  Noetherian, that is, every
$\Gamma$-graded left ideal of $G(A)$ is finitely generated, or equivalently, $G(A)$
satisfies the ascending chain condition for left ideals. Then every finitely
generated $A$-module is left Noetherian, in particular, $A$ is left Noetherian.
\par
(ii) Suppose that $G(A)$ is $\Gamma$-graded left Artinian, that is, $G(A)$
satisfies the descending chain condition for left ideals. Then every finitely
generated $A$-module is left Artinian, in particular, $A$ is left Artinian. \par
(iii) If $G(A)$ is a $\Gamma$-graded simple $K$-algebra, that is, $G(A)$ does not
have nontrivial ideal, then $A$ is a simple $K$-algebra.\par
(iv) Let $M$ be a $\Gamma$-filtered $A$-module with  $\Gamma$-filtration $FM$. If
the Krull dimension (in the sense of Gabriel and Rentschler) of $G(M)$ is
well-defined, then the Krull dimension of $M$ is defined and K.dim$M\le$
K.dim$G(M)$. In particular, this is true for $M=A$ and $G(M)=G(A)$.
\par
(v) Let $M$ be a $\Gamma$-filtered $A$-module with  $\Gamma$-filtration $FM$. If
$G(M)$ is a $\Gamma$-graded simple $G(A)$-module, that is, $G(M)$ does not have
nontrivial $\Gamma$-graded submodule, then $M$ is a simple $A$-module.\par
(vi) If $G(A)$ is semisimple (simple) Artinian, then $A$ is semisimple (simple)
Artinian.\par \vskip 6pt
{\bf Proof} By the foregoing discussion, assertions (i) -- (v) are clear. It
remains to prove the semisimplicity of $A$ in (vi). If $A$ is Artinian, then it is
well-known that the Jacobson radical $J(A)$ of $A$ is nilpotent. As the semisimple
ring $G(A)$ does not contain nilpotent ideal, if we use the $\Gamma$-filtration
$FJ(A)=\{ F_{\gamma}J(A)=J(A)\cap F_{\gamma}A\}_{\gamma\in\Gamma}$ of $J(A)$
induced by $FA$, then $G(J(A))$ is a $\Gamma$-graded ideal of $G(A)$ and hence
$G(J(A))=\{ 0\}$. By Corollary 5.4, $J(A)=\{ 0\}$ as desired. }\QED \v5

\section*{4. Lifting Homological Properties}
In this section we keep the assumption that $\Gamma$ is an ordered {\it monoid}
with the {\it well-ordering} $\prec$, and the identity element $\gamma_0$ of
$\Gamma$  is the {\it smallest} element in $\Gamma$.\par
Let $A$ be a $\Gamma$-filtered $K$-algebra with $\Gamma$-filtration $FA=\{
F_{\gamma}A\}_{\gamma\in\Gamma}$, and let $G(A)$ be the associated $\Gamma$-graded
$K$-algebra of $A$. With notation as before, the aim of this section is to lift
several homological properties of $G(A)$ to $A$, in particular, when $A=R/I$ and
$G(A)=R/\langle\HT (I)\rangle$ as in Theorem 2.1. \par
If $B=\oplus_{\gamma\in\Gamma}B_{\gamma}$ is a $\Gamma$-graded $K$-algebra and
$M=\oplus_{\gamma\in\Gamma}M_{\gamma}$, $N=\oplus_{\gamma\in\Gamma}N_{\gamma}$ are
$\Gamma$-graded $B$-modules, then, by a $\Gamma$-graded $B$-homomorphism from $M$
to $N$ we mean a $B$-homomorphism $\varphi$: $M\r N$ such that $\varphi
(M_{\gamma})\subseteq N_{\gamma}$, $\gamma\in\Gamma$.\v5
We begin with some basics on graded free modules and graded projective modules.\v5
Let $R=\oplus_{\gamma\in\Gamma}R_{\gamma}$ be a $\Gamma$-graded $K$-algebra. A {\it
$\Gamma$-graded free} $R$-module is a free $R$-module $T=\op_{i\in J}Re_i$ on the
basis $\{ e_i\}_{i\in J}$, which is also $\Gamma$-graded such that each $e_i$ is a
homogeneous element, that is, if deg$(e_i)=\gamma_i$, $i\in J$, then
$T=\op_{\gamma\in\Gamma}T_{\gamma}$ with
$$T_{\gamma}=\sum_{i\in J,~w_i\gamma_i=\gamma}R_{w_i}e_i,\quad \gamma\in\Gamma .$$\par
By the definition, to construct a $\Gamma$-graded free $R$-module $T=\oplus_{i\in
J}Re_i$ with the $R$-basis $\{ e_i\}_{i\in J}$, it is sufficient to assign to each
$e_i$ a choosen degree. \v5 Given any $\Gamma$-graded $R$-module
$M=\op_{\gamma\in\Gamma}M_{\gamma}$, $M$ has a generating set $\{ m_i\}_{i\in J}$
consisting of homogeneous elements, i.e., $M=\sum_{i\in J}Rm_i$. Suppose
$d(m_i)=\gamma_i$, $\gamma_i\in\Gamma$, $i\in J$. Then it is easy to see that
$$M_{\gamma}=\sum_{i\in J,~w_i\gamma_i=\gamma}R_{w_i}m_i,\quad\gamma\in\Gamma .$$
Thus, considering the $\Gamma$-graded free $R$-module $T=\oplus_{i\in
J}Re_i=\op_{\gamma\in\Gamma}T_{\gamma}$ with $d(e_i)=\gamma_i$, the map $\varphi$:
$e_i\mapsto m_i$ defines a $\Gamma$-graded $R$-epimorphism $\varphi$: $T\r M$. \v5
If $T$ is a $\Gamma$-graded free $R$-module and $P$ is a $\Gamma$-graded $R$-module,
if there is another $\Gamma$-graded $R$-moduel $Q$ such that $T=P\op Q$ and
$$T_{\gamma}=P_{\gamma}+Q_{\gamma},\quad \gamma\in\Gamma ,$$
then $P$ is called a {\it $\Gamma$-graded projective $R$}-module. \v5 Concerning
graded projective modules, the following result is well-known (e.g., [NVO]).
{\parindent=0pt\v5
{\bf 4.1. Proposition} For a $\Gamma$-graded $R$-module $P$, the following
statements are equivalent.\par
(i) $P$ is a $\Gamma$-graded projective $R$-module.\par
(ii) Given any exact sequence of $\Gamma$-graded $R$-modules and $\Gamma$-graded
$R$-homomorphisms $M\mapright{\psi}{} N\r 0$, if $P\mapright{\alpha}{} N$ is a
$\Gamma$-graded $R$-homomorphism, then there exists a unique $\Gamma$-graded
$R$-homomorphism $P\mapright{\varphi}{} M$ such that the following diagram
commutes:
$$\begin{array}{ccccc}
&&P&&\\
&\scriptstyle{\varphi}\swarrow&\mapdownr{\alpha}&&\\
M&\mapright{}{\psi}&N&\r&0\end{array}\quad \psi\circ\varphi =\alpha$$
(iii) $P$ is projective as an ungraded $R$-module.\par\QED} \v5
Returning to modules over the $\Gamma$-filtered $K$-algebra $A$ with
$\Gamma$-filtration $FA$, we first construct a $\Gamma$-{\it filtered free
$A$-module} $L$ with a $\Gamma$-filtration $FL$ such that its associated
$\Gamma$-graded module $G(L)$ is a $\Gamma$-graded free $G(A)$-module. To this end,
let $L=\op_{i\in J}Ae_i$ be a free $A$-module on the basis $\{ e_i\}_{i\in J}$.
Then, as we did in section 3 (see the proof of Proposition 3.2(ii)), a
$\Gamma$-filtration $FL$ for $L$ can be constructed by using the $A$-basis $\{
e_i\}_{i\in J}$ of $L$ and arbitrarily chosen $\gamma_i\in\Gamma$, $i\in J$, that
is,
$$F_{\gamma}L=\bigoplus_{i\in J}\left (\displaystyle{
\sum_{s\preceq\gamma ,~s_i\gamma_i=s}}F_{s_i}A\right ) e_i,\quad \gamma\in\Gamma
.$$ {\parindent=0pt\v5                                          {\bf 4.2.
Observation} Note that now $\Gamma$ is a monoid with the identity element
$\gamma_0$ which is the smallest element in $\Gamma$. It is not difficult to see
that in the construction of $FL=\{ F_{\gamma}L\}_{\gamma\in\Gamma}$ above, for each
$i\in J$, $e_i\in F_{\gamma_i}L-F_{\gamma_i}^*L$, that is, each $e_i$ is of degree
$\gamma_i$. \v5
{\bf Convention} In what follows, if we say that $L$ is a $\Gamma$-filtered free
$A$-module, then it is certainly the type constructed above. \v5
{\bf 4.3. Proposition} The following statements hold.\par
(i) Let $L=\oplus_{i\in J}Ae_i$ be a $\Gamma$-filtered free $A$-module with
$\Gamma$-filtration $FL$ defined above, then the associated $\Gamma$-graded
$G(A)$-module $G(L)$ of $L$ is a $\Gamma$-graded free $G(A)$-module. More
precisely, we have $G(L)=\oplus_{i\in J}G(A)\sigma
(e_i)=\oplus_{\gamma\in\Gamma}G(L)_{\gamma}$ with
$$G(L)_{\gamma}=\sum_{i\in J,~s_i\gamma_i=\gamma}G(A)_{s_i}\sigma (e_i),\quad \gamma\in\Gamma .$$\par
(ii) If $L'=\oplus_{i\in J}G(A)\eta_i$ is a $\Gamma$-graded free $G(A)$-module with
the $G(A)$-basis $\{\eta_i\}_{i\in J}$ consisting of homogeneous elements, then
there is some $\Gamma$-filtered free $A$-module $L$ such that  $L'\cong G(L)$ as
$\Gamma$-graded $G(A)$-modules.\par
(iii) Let $M$ be a $\Gamma$-filtered $A$-module with  $\Gamma$-filtration $FM=\{
F_{\gamma}M\}_{\gamma\in\Gamma}$. Then there is an exact sequence of
$\Gamma$-filtered $A$-modules and strict $\Gamma$-filtered $A$-homomorphisms
$$0\r N~\mapright{\iota}{}~L~\mapright{\varphi}{}~M\r 0$$
where $L$ is a $\Gamma$-filtered free $A$-module with $\Gamma$-filtration $FL$, $N$
is the kernel of the $\Gamma$-filtered $A$-epimorphism $\varphi$ that has the
$\Gamma$-filtration $FN=\{ F_{\gamma}N=N\cap F_{\gamma}L\}_{\gamma\in\Gamma}$
induced by $FL$, and $\iota$ is the inclusion map.\par
(iv) If $L$ is a $\Gamma$-filtered free $A$-module with $\Gamma$-filtration $F L$,
$N$ is a $\Gamma$-filtered $A$-module with $\Gamma$-filtration $F N$, and
$\varphi$: $G(L)\r G(N)$ is a $\Gamma$-graded $G(A)$-epimorphism, then $\varphi
=G(\psi )$ for some strict $\Gamma$-filtered $A$-epimorphism $\psi$: $L\r N$.
\vskip 6pt \def\GF{G}
{\bf Proof} (i) By the construction of $FL$, Observation 4.2 and the property
$(\sigma )$ of $\sigma$-elements formulated in section 3, the argument is
straightforward.
\par (ii) Suppose $d(\eta_i)=\gamma_i$, $\gamma_i\in\Gamma$, $i\in J$. Then by (i)
we see that the $\Gamma$-filtered free $A$-module $L=\oplus_{i\in J}Ae_i$ with
$d(e_i)=\gamma_i$ satisfies $G(L)\cong L'$.\par (iii) Let $\{\xi_i\}_{i\in J}$ be a
generating set of $M$, that is, $M=\sum_{i\in J}A\xi_i$.  Suppose $\xi_i\in
F_{\gamma_i}M-F_{\gamma_i}^*M$, i.e., $d(\xi_i)=\gamma_i$, $i\in J$. Then the
$\Gamma$-filtered free $A$-module $L=\oplus_{i\in J}Ae_i$ with $d(e_i)=\gamma_i$,
$i\in J$, and the map $\varphi$: $e_i\mapsto \xi_i$ together make the desired exact
sequence.\par (iv) Let $L=\op_{i\in J}Ae_i$ be the $\Gamma$-filtered free
$A$-module with $d(e_i)=\gamma_i$, $i\in J$.  For each $i\in J$, choose $\xi_i\in
F_{\gamma_i}N$ such that $\varphi (\sigma (e_i))= \OV{\xi_i}$, where $\OV{\xi_i}$
is the homogeneous element in $\GF (N)_{\gamma_i}$ represented by $\xi_i$. Then
$\psi$: $L\r N$ can be defined by putting
$$\psi\left (\sum a_ie_i\right )=\sum a_i\xi_i,~\hbox{where}~ \sum a_ie_i\in L.$$
Clearly, $\psi$ is a $\Gamma$-filtered $A$-homomorphism.  Since $G(\psi )$ and
$\varphi$ agree on generators, we have $\GF (\psi)=\varphi$. Hence, by Corollary
3.4, $\psi$ is a strict $\Gamma$-filtered surjection.\QED \v5
{\bf 4.4. Proposition} Let $P$ be a $\Gamma$-filtered $A$-module with
$\Gamma$-filtration $FP=\{ F_{\gamma}P\}_{\gamma\in\Gamma}$. The following
statements hold.\par
(i) If $\GF (P)$ is a projective $\GF (A)$-module, then $P$ is a projective
$A$-module. \par
(ii) If $G(P)$ is a $\Gamma$-graded free $G(A)$-module, then $P$ is a free
$A$-module. \vskip 6pt
{\bf Proof} (i) By Proposition 4.3(iii), there is an exact sequence of
$\Gamma$-filtered $A$-modules and strict $\Gamma$-filtered $A$-homomorphisms
$$0~\mapright{}{}~ N~\mapright{\iota}{}~L~\mapright{\varphi}{}~P~\mapright{}{}~ 0$$
where $L$ is a $\Gamma$-filtered free $A$-module with $\Gamma$-filtration $FL$, $N$
is the kernel of the $\Gamma$-filtered $A$-epimorphism $\varphi$ that has the
$\Gamma$-filtration $FN=\{ F_{\gamma}N=N\cap F_{\gamma}L\}_{\gamma\in\Gamma}$
induced by $FL$, and $\iota$ is the inclusion map. It follows from Proposition 3.3
and Corollary 3.4 that the associated $\Gamma$-graded sequence
$$0~\mapright{}{}~\GF (N)~\mapright{\GF (\iota )}{}~\GF (L)~\mapright{\GF
(\psi)}{}~\GF (P)~\mapright{}{}~0$$ is exact. Since $P$ is a projective
$G(A)$-module, by Proposition 4.1, this sequence splits through $\Gamma$-graded
$G(A)$-homomorphisms. Consequently, $\GF (L)=\GF (P)\op\GF (N)$ with $\GF
(L)_{\gamma}=\GF (P)_{\gamma}\op\GF (N)_{\gamma}$, $\gamma \in\Gamma$, and the
projection of $G(L)$ onto $G(N)$ gives a $\Gamma$-graded $G(A)$-epimorphism $\psi$:
$\GF (L)\r\GF (N)$ such that $\psi\circ\GF (\iota )=1_{\GF (N)}$. Further, by
Proposition 4.3(iv), $\psi =\GF (\beta )$ for some strict $\Gamma$-filtered
$A$-epimorphism $\beta$: $L\r N$. Note that $\GF (\beta )\circ\GF (\iota )=\GF
(\beta\circ\iota )=1_{\GF (N)}$. It follows from Corollary 3.4 that
$\beta\circ\iota$ is an automorphism of $N$. Hence, $L\cong K\op P$. This shows
that $P$ is projective.\par
(ii) Suppose $G(P)=\op_{i\in J}G(A)\sigma (\xi_i)$ with the free $G(A)$-basis
$\{\sigma (\xi_i)\}_{i \in J}$, where each $\xi_i\in
F_{\gamma_i}P-F_{\gamma_i}^*P$, i.e., $d(\xi_i)=\gamma_i\in\Gamma$, $i\in J$.
Then, by Proposition 3.2 (or its proof), $P=\sum_{i\in J}A\xi_i$ with
$$F_{\gamma}P=\displaystyle{\sum_{i\in J}}\left (
\displaystyle{\sum_{s\preceq\gamma ,~s_i\gamma_i=s}}F_{s_i}A \right ) \xi_i,\quad
\gamma\in\Gamma .$$
We claim that $\{\xi_i\}_{i\in J}$ is a free basis for $P$ over $A$. To see this,
construct the $\Gamma$-filtered free $A$-module $L=\op_{i\in J}Ae_i$ with
$\Gamma$-filtration
$$F_{\gamma}L=\bigoplus_{i\in J}\left (\sum_{s\preceq\gamma ,~s_i\gamma_i=s}F_{s_i}A\right )e_i,\quad\gamma\in\Gamma ,$$
as before, such that each $e_i$ has the degree $\gamma_i=d(\xi_i)$. Then we have an
exact sequence of $\Gamma$-filtered $A$-modules and strict $\Gamma$-filtered
$A$-homomorphisms
$$0~\mapright{}{}~N~\mapright{}{}~L~\mapright{\varphi}{}~P~\mapright{}{}~0$$
where $N$ has the $\Gamma$-filtration $FN=\{ F_{\gamma}N=N\cap
F_{\gamma}L\}_{\gamma\in\Gamma}$ induced by $F L$. It follows from Proposition 3.3
that this sequence yields an exact sequence of $\Gamma$-graded $G(A)$-modules and
$\Gamma$-graded $G(A)$-homomorphisms
$$0~\mapright{}{}~G(N)~\mapright{}{}~G(L)~
\mapright{G(\varphi)}{}~G(P)~\mapright{}{}~0$$ However, $G(\varphi )$ is an
isomorphism. Hence $G(N)=\{0\}$ and consequently  $N=\{ 0\}$ by Corollary 3.4. This
proves that $\varphi$ is an isomorphism, or in other words, $P$ is free. \QED \v5
{\bf 4.5. Proposition} Let $M$ be a $\Gamma$-filtered $A$-module with
$\Gamma$-filtration $FM=\{ F_{\gamma}\}_{\gamma\in\Gamma}$, and let
$$0\r N'\r L_n'\r\cdots\r L_0'\r G(M)\r 0\leqno{(1)}$$
be an exact sequence of $\Gamma$-graded $G(A)$-modules and $\Gamma$-graded
$G(A)$-homomorphisms, where the $L_i'$ are $\Gamma$-graded free $G(A)$-modules. The
following statements hold.\par
(i) There exists an exact sequence of $\Gamma$-filtered $A$-modules and strict
$\Gamma$-filtered $A$-homomorphisms
$$0\r N\r L_n\r\cdots\r L_0\r M\r 0\leqno{(2)}$$
in which the $L_i$ are $\Gamma$-filtered free $A$-modules, such that we have the
isomorphism of chain complexes
$$\begin{array}{ccccccccccc}
0\r&N'&\r&L_n'&\r&\cdots&\r&L_0'&\r&G(M)&\r 0\\
&\mapdown{\cong}
&&\mapdown{\cong}&&&&\mapdown{\cong}&&\mapdown{=}&\\
0\r&\GF (N)&\r&\GF (L_n)&\r&\cdots&\r&\GF (L_0)&\r&G(M)&\r 0\end{array}$$
(ii) If $N'$ is a projective $G(A)$-module, then $N$ is a projective $A$-module; If
$N'$ is a $\Gamma$-graded free $G(A)$-module, then $N$ is a free $A$-module.\par
(iii) If all modules in the sequence $(1)$ are finitely generated over $G(A)$, then
all modules in the sequence $(2)$ are finitely generated over $A$. \vskip 6pt
{\bf Proof} (i) By Proposition 4.3, the homomorphism $L_0'\r G(M)$ in sequence
$(1)$ has the form $\GF (\beta )$ for some strict $\Gamma$-filtered surjection
$\beta$: $L_0\r M$, where $L_0$ is a $\Gamma$-filtered free $A$-module such that
$L_0'\cong \GF (L_0)$ as $\Gamma$-graded $G(A)$-modules. Let $N_0=$ Ker$\beta$ and
consider the $\Gamma$-filtration $F N_0=\{ F_{\gamma}N_0=N_0\cap
F_{\gamma}L_0\}_{\gamma\in\Gamma}$ induced by $F L_0$. Then we have the diagram of
$\Gamma$-graded $G(A)$-modules and $\Gamma$-graded $G(A)$-homomorphisms
$$\begin{array}{ccccccccc}
\cdots~\r&L_2'&\r&L_1'&\r&L_0'&\r&G(M)&\r 0\\
         &    &  &    &  &\mapdown{\cong}&&\mapdown{=}&\\
         &0   &\r&\GF (N_0)&\r&\GF (L_0)&\r&G(M)&\r 0\end{array}$$ which has two
         exact rows. Note that the directed square involved in the above diagram
         commutes. It turns out that the homomorphism $L_1'\r L_0'$ factors through
         $\GF (N_0)$, that is, we obtain the diagram
$$\begin{array}{ccccccccc}
\cdots~\r&L_2'&\r&L_1'&\r&L_0'&\r&G(M)&\r 0\\
         &    &  &    \mapdown{}&&\mapdown{\cong}&&\mapdown{=}&\\
         &0   &\r&\GF (N_0)&\r&\GF (L_0)&\r&G(M)&\r 0\\
&&&\downarrow&&&&&\\
&&&0&&&&&\end{array}$$ in which both rows are exact and both directed squares
commute. Starting with $L_1'\r\GF (N_0)\r 0$, the foregoing constructive procedure
can be repeated step by step to yield the desired sequence (2).\par
(ii) and (iii) follow immediately from Proposition 4.3 and Proposition 4.4,
respectively.\QED} \v5
To deal with flat modules over a $\Gamma$-filtered $K$-algebra $A$, we need to
define a $\Gamma$-filtration, respectively a $\Gamma$-gradation, for a tensor
product of two $\Gamma$-filtered $A$-modules, respectively for a tensor product of
two $\Gamma$-graded $G(A)$-modules. \v5 Let $M$ be a $\Gamma$-filtered left
$A$-module with $\Gamma$-filtration $FM$, and let $N$ be a $\Gamma$-filtered right
$A$-module with $\Gamma$ filtration $FN$. Viewing $N\otimes_AM$ as a $\GZ$-module,
we define the $\Gamma$-filtration $F(N\otimes_AM)$ of $N\otimes_AM$ as
$$F_{\gamma}(N\otimes_AM)=\mathbb{Z}\hbox{-span}\left\{ x\otimes y~\Big |~x\in F_vN,~y\in F_wM~\hbox{and}~vw\preceq\gamma\right\},\quad \gamma\in\Gamma .$$
The associated $\Gamma$-graded $\GZ$-module of $N\otimes_AM$ with respect to
$F(N\otimes_AM)$ is then defined as
$G(N\otimes_AM)=\oplus_{\gamma\in\Gamma}G(N\otimes_AM)_{\gamma}$ with
$$G(N\otimes_AM)_{\gamma}=F_{\gamma}(N\otimes_AM)/F_{\gamma}^*(N\otimes_AM),\quad\gamma\in\Gamma ,$$
where $F_{\gamma}^*(N\otimes_AM)=\cup_{\gamma'\prec\gamma}F_{\gamma'}(N\otimes M)$.
\par Let $P$ be a $\Gamma$-graded left $G(A)$-module, and let $Q$ be a
$\Gamma$-graded right $G(A)$-module. Viewing $Q\otimes_{G(A)}P$ as a $\GZ$-module,
we define the $\Gamma$-gradation of $Q\otimes_{G(A)}P$ as
$$(Q\otimes_{G(A)}P)_{\gamma}=\mathbb{Z}\hbox{-span}\left\{ z\otimes t~\Big |~z\in Q_v,~t\in P_w~\hbox{and}~vw=\gamma\right\},\quad \gamma\in\Gamma .$$
\v5 {\parindent=0pt\par
{\bf 4.6. Lemma} Let $M$ be a $\Gamma$-filtered left $A$-module with
$\Gamma$-filtration $FM$, and let $N$ be a $\Gamma$-filtered right $A$-module with
$\Gamma$ filtration $FN$. With the definition made above, the following statements
hold.\par (i) For $\OV x_v\in G(N)_v$ represented by $x\in F_vN$, and $\OV y_w\in
G(M)_w$ represented by $y\in F_wM$, the mapping defined by
$$\begin{array}{cccc} \varphi (M,N):&G(N)\otimes_{G(A)}G(M)&\mapright{}{}&G(N\otimes_AM)\\
&\OV x_v\otimes \OV y_w&\mapsto&(\overline{x\otimes y})_{vw}\end{array}$$ is an
epimorphism of $\Gamma$-graded $\mathbb{Z}$-modules.\par (ii) The canonical
$A$-isomorphisms
$$A\otimes_AM~\mapright{\cong}{}~M~\hbox{and}~N\otimes_AA~\mapright{\cong}{}~N$$
are strict $\Gamma$-filtered $A$-isomorphisms. \par (iii) The strict
$\Gamma$-filtered $A$-isomorphisms in (ii) induce $\Gamma$-graded
$G(A)$-isomorphisms
$$G(A\otimes_AM)~\mapright{\cong}{}~G(M)~\hbox{and}~G(N\otimes_AA)~\mapright{\cong}{}~G(N).$$
(iv) The canonical $G(A)$-isomorphisms
$$G(A)~\otimes_{G(A)}G(M)~\mapright{\cong}{}~G(M)~\hbox{and}~G(N)\otimes_{G(A)}G(A)~\mapright{\cong}{}~G(N)$$
are $\Gamma$-graded $G(A)$-isomorphisms. \vskip 6pt {\bf Proof} Verification is
straightforward.\QED \v5
{\bf 4.7. Proposition} Let $M$ be a $\Gamma$-filtered left $A$-module with
$\Gamma$-filtration $FM$. If $G(M)$ is a flat $\Gamma$-graded $G(A)$-module, then
$M$ is a flat $A$-module. \vskip 6pt
{\bf Proof} Let $J$ be a right ideal of $A$ and $FJ=\{ F_{\gamma}J=J\cap
F_{\gamma}A\}_{\gamma\in\Gamma}$ the $\Gamma$-filtration of $J$ induced by $FA$.
Consider the inclusion map $\iota$: $J\hookrightarrow A$. Then the strict exactness
of the $\Gamma$-filtered sequence
$$0~\mapright{}{}~J~\mapright{\iota}{}~A$$
yields the exact $\Gamma$-graded sequence
$$0~\mapright{}{}~G(J)~\mapright{G(\iota )}{}~G(A)$$
Furthermore, it follows from the flatness of $G(M)$ and Lemma 4.6 that we have the
following commutative diagram of $\Gamma$-graded $\mathbb{Z}$-modules and
$\Gamma$-graded $\mathbb{Z}$-homomorphisms:
$$\begin{array}{cccccc}
0\r&G(J)\otimes_{G(A)}G(M)&\mapright{G(\iota )\otimes1_{G(M)}}{}&G(A)\otimes_{G(A)}G(M)&&\\
&&&&\searrow\scriptstyle{\cong}&\\
&\mapdown{\varphi (M,J)}&&\mapdown{\varphi (M,A)}&&G(M)\\
&&&&\swarrow\scriptstyle{\cong}&\\
&G(J\otimes_AM)&\mapright{}{G(\iota\otimes1_M)}&G(A\otimes_AM)&&\end{array}$$
As $\varphi (M,A)$ is an isomorphism and $G(\ell )\otimes 1_{G(M)}$ is a
monomorphism, it turns out that $\varphi (M,J)$ is an isomorphism. So
$G(\iota\otimes 1_M)$ must be a monomorphism. By previous Corollary 3.4, we
conclude that $\iota\otimes 1_M$ is a strict $\Gamma$-filtered monomorphism. This
proves the flatness of $M$.\QED} \v5
Let $A$ be a $\Gamma$-filtered $K$-algebra with $\Gamma$filtration $FA$. Noticing
every $A$-module can be endowed with a $\Gamma$-filtration (section 1 Example (2)),
Proposition 4.5 and Proposition 4.7 enable us to reach the main results of this
section. In the text below we write p.dim to denote the projective dimension of a
module, gl.dim to denote  the homological global dimension of a ring, and gl.w.dim
to denote the global weak dimension of a ring, respectively; and moreover, we write
w.dim for the weak dimension of a module. {\parindent=0pt\v5
{\bf 4.8. Theorem} Let $A$ be a $\Gamma$-filtered $K$-algebra with
$\Gamma$-filtration $FA$, and let $G(A)$ be the associated $\Gamma$-graded
$K$-algebra of $A$. The following statements hold, especially when $A=R/I$ and
$G(A)=R/\langle\HT (I)\rangle$ as in Theorem 2.1. \par
(i) Let $M$ be an $A$-module with $\Gamma$-filtration $FM$. Then p.dim$_AM\le$
p.dim$_{G(A)}G(M)$. In particular, if $G(M)$ has a (finite or infinite
$\Gamma$-graded) free resolution, then $M$ has a (finite or infinite) free
resolution.\par
(ii) gl.dim$A\le$ gl.dom$G(A)$.\par
(iii) If $G(A)$ is left hereditary, then $A$ is left hereditary.\par
(iv) Let $M$ be an $A$-module with $\Gamma$-filtration $FM$. Then w.dim$_AM\le$
w.dim$_{G(A)}G(M)$.\par
(v) $\hbox{gl.w.dim}A\le \hbox{gl.w.dim}G(A)$.\par
(vi) If $G(A)$ is a Von Neuman regular ring then so is $A$. \vskip 6pt
{\bf Proof} (i) and (ii) are immediate consequences of Proposition 4.4 and
Proposition 4.5.\par
(iii) If $G(A)$ is left hereditary, then every left ideal of $G(A)$ is a projective
$G(A)$-module. Let $L$ be a left ideal of $A$ and $FL$ the $\Gamma$-filtration of
$L$ induced by $FA$. Using the inclusion map $L\hookrightarrow A$ (note that this
is a strict $\Gamma$-filtered $A$-homomorphism), we may view $G(L)$ as a
$\Gamma$-graded left ideal of $G(A)$ by Corollary 3.4. Thus, $G(L)$ is a projective
$G(A)$-module, and it follows from Proposition 4.4 that $L$ is a projective
$A$-module. Therefore, $A$ is left hereditary.\par
(iv) Note that any exact sequence
$$0\r N\r L_n \r\cdots\r L_1\r L_0\r M\r 0$$
consisting of $\Gamma$-filtered free $A$-modules $L_i$ and strict $\Gamma$-filtered
$A$-homomorphisms yields an exact sequence
$$0\r G(N)\r G(L_n) \r\cdots\r G(L_1)\r G(L_0)\r G(M)\r 0$$
consisting of $\Gamma$-graded free $G(A)$-modules $G(L_i)$ and $\Gamma$-graded
$G(A)$-homomorphisms, where $N$ has the $\Gamma$-filtration $FN$ induced by $FL_n$.
This assertion is an immediate consequences of Proposition 4.7.\par
(v) and (vi) follow from (iv).\QED} \v5

\section*{5. With Gr\"obner Bases: $G^{\B}(R/I)\cong  R/\langle\LM (\G )\rangle$ \&
$G^{\NZ}(R/I)\cong R/\langle\HT (\G )\rangle$}                                  By
using Gr\"obner bases in a computational setting, the aim of this section is to
decode the defining relations of the associated graded $K$-algebra of the
$K$-algebra $R/I$ in Theorem 2.1 with respect to the $\B$-filtration and the
$\NZ$-filtration of $R/I$, respectively. \v5
Let $R=K[a_1,...,a_n]$ be a finitely generated $K$-algebra over a field $K$, where
$R$ has a $K$-basis $\B$  consisting of {\it monomials} of the form
$$u=a_{i_1}\cdots a_{i_s},\quad a_{i_j}\in\{ a_1,...,a_n\} ,~s\in\NZ ,~s\ge 1.$$
Suppose that $\B$ is a {\it skew multiplicative $K$-basis} of $R$ in the sense that
$$u,~v\in\B~\hbox{implies}~\left\{\begin{array}{l} u\cdot v=\lambda w~\hbox{for some}
~\lambda\in K^*,~w\in\B ,\\
\hbox{or}~u\cdot v=0.\end{array}\right.\leqno{(\hbox{sm})}$$
The reason that we use the word ``skew'' here is that free algebras, commutative
polynomial algebras, the coordinate rings of quantum affine $K$-spaces, and path
algebras defined by finite directed graphs, all are involved as the most important
practical examples supporting our text. \par
Let $\prec$ be a total ordering on $\B$. If we adopt the commonly used terminology
in computational algebra, then for $f\in R$, say
$$f=\sum^s_{i=1}\lambda_iu_i,\quad \lambda_i\in K^*,~u_i\in\B,~u_1\prec u_2\prec\cdots\prec u_s ,$$
the {\it leading monomial} of $f$, denoted $\LM (f)$, is defined as $\LM (f)=u_s$;
the {\it leading coefficient} of $f$, denoted $\LC (f)$, is defined as $\LC (f)
=\lambda_s$; and the {\it leading term} of $f$, denoted $\LT (f)$, is defined as
$\LT (f)=\LC (f)\LM (f)=\lambda_su_s$.  Thus, for a subset $S$ of $R$, the set of
leading monomials of $S$ is defined as  $\LM (S)=\{ \LM (f)~|~f\in S\}$. \par
Under the assumption (sm) on $\B$, recall that a {\it monomial ordering} on $R$ is
a {\it well-ordering} $\prec$ on $\B$ satisfying the following
conditions:{\parindent=1.5truecm\vskip 6pt
\re{{\bf (Mo1)}} If $u\prec v$, then $\LM (uw)\prec \LM (vw)$ if both $uw\ne 0$ and
$vw\ne 0$.
\re{{\bf (Mo2)}} If $u\prec v$, then $\LM (su)\prec \LM (sv)$ if both $su\ne 0$ and
$sv\ne 0$.
\re{{\bf (Mo3)}} If $uw=\lambda v$, then $v\succ u$ and $v\succ
w$.}{\parindent=0pt\par
Besides, if $1\in\B$, then it is required that $1\prec u$ for all $u\in\B-\{ 1\}$,
and moreover, $v,u,w\ne 1$ in the axiom (Mo3).}\vskip 6pt
If $\prec$ is a monomial ordering on $R$, then, by mimicking (e.g., [Gr]), $R$
holds a Gr\"obner basis theory, that is, theoretically every  ideal $I$ of $R$ has
a (finite or infinite) Gr\"obner basis $\G$ in the sense that
$$\langle\LM (I)\rangle =\langle\LM (\G )\rangle .$$
{\parindent=0pt\v5
{\bf $\B$-filtered case}\par
Let $R$ be as fixed above. In this part we assume that $R$ {\it does not have
divisors of zero}, $1\in\B$ (thus, path algebras are excluded), and le $\prec$ be a
monomial ordering on $R$. Hence $R$ holds a Gr\"obner basis theory with respect to
$(\B ,\prec )$.}
\par
Note that $R$ is $\B$-graded, namely, $R=\oplus_{u\in\B}R_u$ with $R_u=Ku$. In this
case, we see that for $f\in R$, the head term $\HT(f)$ of $f$ defined in section 2
is the same as the leading term $\LT (f)$ of $f$ defined above, that is,
$$\HT(f)=\lambda_su_s=\LT (f)=\LC (f)\LM (f),$$
It turns out that if $I$ is an ideal of $R$, then
$$\langle\HT (I)\rangle =\langle\LM (I)\rangle ,$$
where the latter is usually called the {\it initial monomial ideal} of $I$, and
consequently $R/\langle\LM (I)\rangle$ is called the {\it associated monomial
algebra} of the algebra $R/I$.  \par
Since $R$ has no divisors of zero, by section 1,  $R$ is $\B$-filtered by the
$\B$-grading filtration $F^{\B}R=\{ F^{\B}_uR\}_{u\in\B}$, where
$$F^{\B}_uR=\oplus_{v\preceq u}R_v,\quad u\in\B .$$
Now, let $I$ be an ideal of $R$ and $A=R/I$, the quotient algebra of $R$ defined by
$I$. Then $A$ has the $\B$-filtration $F^{\B}A=\{ F^{\B}_uA\}_{u\in\B}$ induced by
$FR$, that is,
$$F^{\B}_uA=(F^{\B}_uR+I)/I,\quad u\in\B , $$
which defines the associated $\B$-graded $K$-algebra
$G^{\B}(A)=\oplus_{u\in\B}G^{\B}(A)_u$  with $G^{\B}(A)_u=F^{\B}_uA/F^{\B
*}_{~~u}A$ (see section 1).  {\parindent=0pt \v5
{\bf 5.1. Theorem} With notation as fixed above, let $\G$ be a generating set of
the  ideal $I$. The following statements are equivalent. \par
(i) $\G$ is a Gr\"obner basis for $I$ with respect to the given monomial ordering
$\prec$ on $\B$.\par
(ii) $\langle\LM (I)\rangle =\langle\LM (\G )\rangle$.\par
(iii) $G^{\B}(A)\cong R/\langle\LM (I)\rangle =R/\langle\LM (\G )\rangle$.\vskip
6pt     {\bf Proof} Note that $\langle\LM (\G )\rangle\subseteq \langle\LM
(I)\rangle$. This follows immediately from the definition of a Gr\"obner basis in
$R$ and Theorem 2.1. \QED \v5
{\bf Remark} In view of Theorem 2.1, Theorem 3.1(i) and Theorem 5.1, actually, a
richer  Gr\"obner basis theory in both commutative and noncommutative cases may be
introduced by solving the isomorphic problem
$$G^{\B}(A)~\mapright{\cong}{\hbox{?}}~R/\langle\LM (F)\rangle$$
for a given generating set $F$ of the ideal $I$. On this aspect, a systematic
clarification has been done in [Li3]. \v5
{\bf $\NZ$-graded case}\par
In this part we {\it allow the case} that $R$ has divisors of zero, for instance,
$R$ is a path algebra defined by a finite directed graph. }\par
By the choice of the $K$-basis $\B$, $R$  is also $\NZ$-graded by the natural
$\NZ$-gradation $\{ R_p\}_{p\in\NZ}$ defined by lengths  of elements in $\B$, that
is, $R=\oplus_{p\in\NZ}R_p$ with
$$R_p=K\hbox{-span}\left\{ u=a_{i_1}^{\alpha_1}\cdots a_{i_s}^{\alpha_s}\in\B~\Big |~
\alpha_1+\cdots +\alpha_s=p\right\} ,~p\in\NZ .$$                             Also
recall that if $f\in R$, $f=r_p+r_{p-1}+\cdots +r_0$ with $r_i\in R_i$ and $r_p\ne
0$, then the head term of $f$ is defined as $\HT (f)=r_p$, and we say that $f$ is
of degree $p$ in $R$, denoted $d(f)=p$. For a subset $S\subset R$, we write
$$\HT (S)=\left\{ \HT (f)~\Big |~f\in S\right\} .$$\par
Further, let $\prec$ be a well-ordering on $\B$. If the ordering $\prec_{gr}$
defined for $u$, $v\in\B$ by the rule
$$\begin{array}{rcl} u\prec_{gr} v&\Leftrightarrow& d(u)<d(v)\\
&~~&\hbox{or}~d(u)=d(v)~\hbox{and}~u\prec v\end{array}$$
is a monomial ordering on $\B$ in the sense of the foregoing (Mo1) -- (Mo3), then
we call $\prec_{gr}$ a {\it graded monomial ordering} on $\B$.
 {\parindent=0pt\v5
{\bf 5.2. Theorem} (A generalization of [Li1] CH.III Theorem 3.7) Let $\prec_{gr}$
be a graded monomial ordering on $\B$, and let $I$ be an ideal of $R$. Put
$J=\langle\HT (I)\rangle$. The following two statements hold.\par
(i) $\LM (J)=\LM (I)$.\par
(ii) Let $\G$ be a generating set of $I$. Then $\G$ is a Gr\"obner basis of $I$
with respect to $(\B ,\prec_{gr})$ if and only if  $\HT (\G )$ is a Gr\"obner basis
for the $\NZ$-graded ideal $J$ of $R$ with respect to $(\B ,\prec_{gr})$. \vskip
6pt
{\bf Proof} (i) First, note that $\prec_{gr}$ is a graded monomial ordering on
$\B$. For $f\in R$, we have
$$\LM (f)=\LM (\HT (f)),\leqno{(*)}$$
and this turns out $\LM (I)=\LM (\HT (I ))$.  Hence $\LM (I)\subset\LM (J)$. It
remains to prove the inverse inclusion. Since $J$ is an $\NZ$-graded ideal of $R$,
noticing the formula $(*)$ above, we need only to consider the leading monomials of
homogeneous elements. Let $F\in J$ be a homogeneous element of degree $p$. Then
$F=\sum_iG_i\HT (f_i)H_i$, where $G_i$, $H_i$ are homogeneous elements of $R$ and
$f_i\in I$, such that $d(G_i)+d(f_i)+d(H_i)=p$ whenever $G_i\HT (f_i)H_i\ne 0$.
Write $f_i=\HT (f_i)+f_i'$ such that $d(f_i')<d(f_i)$. Then
$$\sum_iG_if_iH_i=F+\sum_iG_if_i'H_i,$$
in which $d(\sum_iG_if_i'H_i)<p=d(F)$. Hence $\LM (F)=\LM (\sum_iG_if_iH_i)\in\LM
(I)$. This shows that $\LM (J)\subset \LM (I)$, and consequently, the desired
equality follows.\par
(ii)  Note that the above formula $(*)$ yields $\langle\LM (\G )\rangle =\langle\LM
(\HT (\G ))$. By the equality in (i) we have
$$ \langle\LM (I)\rangle =\langle\LM (\G )\rangle ~\hbox{if and only if}~
\langle\LM (J)\rangle =\langle\LM (\HT (\G ))\rangle .$$
Hence the equivalence follows. \QED}\v5
Let $I$ be an ideal of $R$ and $A=R/I$. Then it follows from section 1 that the
$\NZ$-grading filtration $F^{\NZ}R=\{ F^{\NZ}_pR\}_{p\in\NZ}$ of $R$ with
$F^{\NZ}_pR=\oplus_{i\le p}R_i$ induces the natural $\NZ$-filtration $F^{\NZ}A=\{
F^{\NZ}_pA\}_{p\in\NZ}$ of $A$ with
$$F^{\NZ}_pA=(F^{\NZ}_pR+I)/I,\quad p\in\NZ,$$
that defines the associated $\NZ$-graded $K$-algebra $G^{\NZ}(A)$ of $A$, namely,
$$G^{\NZ}(A)=\bigoplus_{p\in\NZ}G^{\NZ}(A)_p~\hbox{with}~G^{\NZ}(A)_p=F^{\NZ}_pA/F^{\NZ}_{p-1}A.$$
{\parindent=0pt\par
{\bf 5.3. Theorem}  (A generalization of [Li1] CH.III Theorem 3.6, [Li2] Theorem
2.1) Let $\prec_{gr}$ be a graded monomial ordering on $\B$ and let $I$ be an ideal
of $R$. Consider the algebra $A=R/I$ with the natural $\NZ$-filtration $F^{\NZ}A=\{
F^{\NZ}_pA\}_{p\in\NZ}$, and let $G^{\NZ}(A)$ be the associated $\NZ$-graded
algebra of $A$. Suppose $\G$ is a Gr\"obner basis of $I$ with respect to $(\B
,\prec_{gr})$. Then $\langle\HT (I)\rangle =\langle\HT (\G )\rangle $ and hence
$$G^{\NZ}(A)\cong R/\langle\HT (I)\rangle =R/\langle\HT (\G )\rangle .$$
{\bf Proof}  We prove the equality $\langle\HT (I)\rangle =\langle\HT (\G )\rangle
$ by showing that $\G$ satisfies the condition of Proposition 2.2(i). Let $f\in I$.
As $\G$ is a Gr\"obner basis for $I$, $f$ has a presentation
$$f=\sum_j\lambda_ju_jg_jv_j,\quad \lambda_j\in K,~u_j,v_j\in\B,~g_j\in\G ,$$
satisfying $\LM (u_j\LM (g_j)v_j)\preceq_{gr}\LM (f)$ for all $u_jg_jv_j\ne 0$ this
is a result by division by $\G$). Since $\prec_{gr}$ is a graded monomial ordering
on $\B$, for any $h\in R$ we have $\LM (h)=\LM (\HT (h))$ and $d(\LM (h))=d(\HT
(h))=d(h)$. Thus, the presentation of $f$ yields
$$d(u_j)+d(g_j)+d(v_j)\le d(f)~\hbox{for all}~u_jg_jv_j\ne 0,$$
as desired. Now it follows from Theorem 2.1 that $G^{\NZ}(A)\cong R/\langle\HT
(I)\rangle =R/\langle\HT (\G )\rangle $.\QED}  \v5
By Theorem 5.2(i), in the case that a graded monomial ordering $\prec_{gr}$ is
used, we see that if $R/I$ has the $\B$-filtration induced by the $\B$-grading
filtration $FR$ of $R$, then
$$G^{\B}(R/I)\cong R/\langle\LM (I)\rangle\cong G^{\B}(R/\langle\HT (I)\rangle =G^{\B}(G^{\NZ}(R/I)) .$$
Combining the classical lifting technics for $\NZ$-filtered algebras, we now
summarize the whole lifting strategy of sections 2 -- 4 with respect to both
$\B$-filtration and $\NZ$-filtration in the following diagram:\par
$$\begin{diagram} &&R/I&&\\
&\NE^{\textsf{lifting}}&&\NW^{\textsf{lifitng~(classical)}}&\\
G^{\B}(R/I)&&&&G^{\NZ}(R/I)=\FRAC{R}{\langle\HT (I)\rangle}=\FRAC{R}{\langle\HT (\G )\rangle}\\
\uTo^{\cong}&&&&\uTo_{\textsf{~lifting}}\\
\FRAC{R}{\langle\LM (\G )\rangle}=\FRAC{R}{\langle\LM
(I)\rangle}&&\rTo_{\cong}&&G^{\B}\left (\FRAC{R}{\langle\HT
(I)\rangle}\right )\\
\end{diagram}$$
\v5
\section*{6. The first application}\par
As the first application of Theorem 5.1 and Theorem 5.3, we now can mention a
result that generalizes several well-known facts.\par
Let $R=K[a_1,...,a_n]$ and the skew multiplicative $K$-basis $\B$ of $R$ be as in
section 5. Suppose that $\prec$ is a monomial ordering on $\B$. If $I$ is an ideal
of $R$, then by the fundamental decomposition theorem for the vector space $R$, we
have
$$R=I\oplus K\hbox{-span}(\B -\LM (I))=\langle\LM (I)\rangle\oplus K\hbox{-span}(\B -\LM (I)) .$$
{\parindent=0pt
{\bf 6.1. Corollary} With assumption as made above, Let $A=R/I$, and let $F^{\B}A$
be the $\B$-filtration of $A$ (if it exists) and $F^{\NZ}A$ the natural
$\NZ$-filtration of $A$. Write $G^{\B}(A)$ and $G^{\NZ}(A)$ for the associated
graded algebras of $A$ with respect to $F^{\B}A$ and $F^{\NZ}A$, respectively. The
following statements hold.\par
(i) The image of the set $\B -\LM (I)$ in  $A=R/I$, $G^{\B}(A)= R/\langle\LM
(I)\rangle$, and $G^{\NZ}(A)= R/\langle\HT (I)\rangle$, respectively, serves to
give a $K$-basis for each algebra listed.\par
(ii) $A=R/I$ is finite dimensional over $K$ if and only if $G^{\B}(A)=R/\langle\LM
(I)\rangle$ is finite dimensional over $K$ if and only if $G^{\NZ}(A)=R/\langle\HT
(I)\rangle$ is finite dimensional over $K$, and in this case we have
$$\hbox{dim}_KA=\hbox{dim}_KG^{\B}(A)=\hbox{dim}_KG^{\NZ}(A)=|\B -\LM (I)|.$$\par
(iii) Consider the natural $\NZ$-filtration for $A=R/I$, $G^{\B}(A)= R/\langle\LM
(I)\rangle$, and $G^{\NZ}(A)= R/\langle\HT (I)\rangle$,respectively. Then all three
$K$-algebras have the same Hilbert function, and hence, they have the same growth,
or equivalently, they have the same Gelfand-Kirillov dimension. \par
(iv) If $I$ is an $\NZ$-graded ideal of $R$, then the $\NZ$-graded $K$-algebra
$A=R/I$ and the $\B$-$\NZ$-graded $K$-algebra $G^{\B}(A)=R/\langle\LM (I)\rangle$
have the same Hilbert series. Hence in (iii) above, $G^{\B}(A)= R/\langle\LM
(I)\rangle$ and $G^{\NZ}(A)= R/\langle\HT (I)\rangle$ have the same Hilbert series.
(v) If $\G$ is a Gr\"obner basis of $I$, then the set $\B -\LM (I)$, the
Gelfand-Kirillov dimension, and the Hilbert series of the respective algebra
considered in (i) -- (iv) above  may be obtained algorithmically.\par}
\par\QED\v5

\section*{7. Realization via Gr\"obner Bases and Ufnarovski Graphs}
Thanks to [[An1], [An2], [G-IL], [G-I1], [G-I2], [Nor], [Uf1], and [Uf2], in this
section we indicate how to realize some of the foregoing lifting properties by
virtue of Gr\"obner bases and the associated Ufnarovski (chain) graphs. To be
concrete, we focus on a free $K$-algebra $R=\KS$ with $X=\{ X_1,...,X_n\}$ and the
data $(\B ,\prec_{gr})$, where $\B$ is the standard $K$-basis and $\prec_{gr}$ is
some graded monomial ordering on $\B$. All notations are retained as before.\v5
Let $\Omega=\{ u_1,...,u_s\}$ be a finite {\it reduced} subset of $\B$ in the sense
that $u_i$ and $u_j$ are not divisible each other if $i\ne j$. If $u_i\in\Omega$,
then, as before we write $d(u_i)$ for the degree (length) of $u_i$ with respect to
the $\NZ$-gradation of $R$. Put
$$\ell =\max\left\{ d(u_i)~\Big |~u_i\in\Omega\right\} .$$
 Then the {\it Ufnarovski graph}  of $\Omega$ (in the
sense of [Uf1]), denoted $\Gamma (\Omega)$, is defined as a directed graph, in
which the set of vertices $V$ is given by
$$V=\left\{ v_i~\Big |~v_i\in\B -\langle\Omega\rangle,~d(v_i)=\ell -1\right\} ,$$
and the set of edges $E$ contains the edge $v_i\r v_j$ if and only if there exist
$X_i$, $X_j\in X$ such that $v_iX_i=X_jv_j\in\B -\langle\Omega\rangle$. The
essential link between this graph and the monomial $K$-algebra
$R/\langle\Omega\rangle$ is that there is a bijective correspondence between the
set of monomials of degree $\ge \ell -1$ in $\B -\langle\Omega\rangle$ and paths in
the graph.\par The first effective application of Ufnarovski graph was made to
determine the growth of a finitely presented $K$-algebra, for instance,
{\parindent=.5truecm\re{} \re{$\bullet$} ([Uf1] 1982) the growth of
$R/\langle\Omega\rangle$ is exponential if and only if there are two different
cycles in the graph $\Gamma (\Omega )$ with a common vertex. Otherwise,
$R/\langle\Omega\rangle$ has polynomial growth of degree $d$, where $d$ is the
maximal possible number of different cycles in $\Gamma (\Omega )$ through which one
path can pass. Hence, as a $K$-vector space, dim$_K(R/\langle\Omega\rangle )
<\infty$ if and only if the graph $\Gamma (\Omega )$ does not contain any
cycle.\re{}} {\parindent=0pt\par In view of section 6, if $A=R/I$ is a finitely
presented $K$-algebra (i.e., $I$ is a finitely generated ideal of $R$), and if $\G$
is a finite reduced Gr\"obner basis for the ideal $I$ (it is known that if $I$ has
a Gr\"obner basis, then a reduced Gr\"obner basis can always be obtained by using
the division algorithm), then the same statement as above can be mentioned for $A$
and its associated $\NZ$-graded algebra $G^{\NZ}(A)=R/\langle\HT (\G )\rangle$ by
means of the Ufnarovski graph $\Gamma (\LM (\G ))$. }\v5 Let $I$ be an arbitrary 
ideal of $R$ and $\HT (I)$ the set of head terms of $I$ with respect to the natural 
$\NZ$-gradation of $R$. Below we realize other lifting properties for algebras 
$R/I$ and $R/\langle\HT (I)\rangle$ ($\cong G^{\NZ}(R/I)$).
\par
 {\parindent=0pt\v5
{\bf Noetherianity}\par                                                        
First of all, by Theorem 3.5, the following proposition is clear. \v5                                                                      
{\bf 7.1. Proposition} Let  $J$ be a monomial ideal of $R$ such that $J\subseteq 
\langle\LM (I)\rangle$. If the monomial algebra $R/J$ is left (right) Noetherian 
then the algebra $R/I$ and the algebra $R/\langle\HT (I)\rangle$ are left (right) 
Noetherian.
\par\QED \v5
{\bf 7.2. Theorem} With notation as fixed above, let $R/J$ be a finitely presented
monomial algebra defined by the monomial ideal $J=\langle\Omega\rangle$ with
$\Omega =\{ u_1,...,u_s\}\subset\B$ a reduced subset. Suppose that $J\subseteq
\langle\LM (I)\rangle$ (for instance, $\Omega\subseteq\LM (I)$).  If there is no
edge entering (leaving) any cycle of the graph $\Gamma (\Omega )$, then the algebra
$A=R/I$ and $R/\langle\HT (I)\rangle$ are left (right) Noetherian.\vskip 6pt {\bf
Proof} By ([Uf2], [Nor]), the finitely presented monomial algebra $R/J$ is left
(right) Noetherian if and only if there is no edge entering (leaving) any cycle of
the graph $\Gamma (\Omega )$. So the theorem follows from Proposition 7.1.\par\QED}
\v5
As a small example let us look at the monomial algebra $S=K\langle X,Y\rangle
/\langle X^2, YX\rangle$. Put $\Omega =\{ X^2,YX\}$. It is easy to see that the
graph $\Gamma (\Omega )$ is of the form
$$\begin{array}{ccl} X&&Y\\ \bullet&\mapright{}{}&\bullet\circlearrowleft\end{array}$$
and there is no edge leaving the only cycle of $\Gamma (\Omega )$. Hence $S$ is
right Noetherian. It follows from Theorem 7.3 that  any algebra $A=K\langle
X,Y\rangle /I$ with $X^2$, $YX\in \langle\LM (I)\rangle$ and the algebra $K\langle
X,Y\rangle /\langle\HT (I)\rangle$ are right Noetherian, for instance, if we set
$Y\prec_{gr} X$ and take $I=\langle X^2+aY^2+b,YX+cX+d,h(X,Y)\rangle$, where
$a,b,c,d\in K$, $h(X,Y)\in K\langle X,Y\rangle$. {\parindent=0pt\v5             
{\bf Remark} Recall that an algebra  is called {\it weak Noetherian} if it 
satisfies the ascending chain condition for ideals. If $\Omega\subset\B$ is a 
finite reduced subset, then it was proved in [Nor] that the algebra 
$R/\langle\Omega\rangle$ is weak Noetherian if and only if the Ufnarofski graph 
$\Gamma (\Omega )$ does not contain any cycle with edges both entering and leaving 
it. In a similar way one may also get an analogue of Theorem 7.2 on the weak 
Noetherianity of $R/I$ and $R/\langle\HT (I)\rangle$. \v5                                                  
{\bf Semisimplicity, primeness and semiprimeness}\par                            
Let $\G=\{ g_1,...,g_s\}$ be a reduced Gr\"obner basis of $I$ with respect to $(\B 
,\prec_{gr})$, $\LM (\G )$ the set of leading monomials of $\G$, and let $\ell 
=\max\{ d(u)~\Big |~u\in\LM (\G )\}$. Recall from [Uf1] and [G-I2] that a vertex in 
the Ufnarovski graph $\Gamma (\LM (\G ))$ is called {\it cyclic} if it belongs to a 
cyclic route of $\Gamma (\LM (\G ))$. For each monomial $v=x_{i_1}x_{i_2}\cdots 
x_{i_s}$ with $s>\ell -1$, there is a unique route of $\Gamma (\LM (\G ))$ which is 
defined by
$$R(v)=v_0\r v_1\r\cdots \r v_d,$$
where $d=s-\ell$ and $v_j=x_{i_{j+1}}x_{i_{j+2}}\cdots x_{i_{j+\ell}}$, $0\le j\le 
d$. A monomial $v\in\B -\langle\LM (\G )\rangle$, $v\ne 1$, is called {\it cyclic} 
if $d(v)\le \ell -1$ and $v$ is a right-hand segment of a cyclic vertex in $\Gamma 
(\LM (\G ))$, or, if $d(v)>\ell -1$ and the route $R(v)$ is a subroute of a cyclic 
route. \v5                                                                      
{\bf 7.3. Theorem} With convention made above, the following statements hold.\par 
(i) If any $v\in\B -\langle\LM (\G )\rangle$ with $1\le d(v)\le \ell$ is cyclic, 
then $R/I$ and $R/\langle\HT (I)\rangle$ are semiprime; If furthermore $I$ is not 
an $\NZ$-graded ideal of $R$ and $R/I$ is artinian (for example 
dim$_K(R/I)<\infty$), then $R/I$ is semisimple artinian.\par                    
(ii) If $\Gamma (\LM (\G ))$ satisfies\par (a) any $v\in\B -\langle\LM (\G 
)\rangle$ with $d(v)<\ell -1$ is a right-hand segment of a vertex of $\Gamma (\LM 
(\G ))$, and 
\par (b) for any two vertices $u$ and $v$ of $\Gamma (\LM (\G ))$, there exists a 
route from $u$ to $v$,\par                                                      
then $R/I$ and $R/\langle\HT (I)\rangle$ are prime rings; If furthermore $I$ is not 
an $\NZ$-graded ideal of $R$ and $R/I$ is artinian , then $R/I$ is semisimple 
artinian. \vskip 6pt                                                            
{\bf Proof} This follows from Theorem 3.1 and ([G-I1] Theorem 2.21, 2.27 and 
2.28).\QED} \v5                                                               
Recall that for any $p\in\NZ$, $\B_p=\{ w\in\B~|~d(w)=p\}$ is a finite set. Also 
note that for a monomial $w\in\B$, the property that $w\in\B -\langle\LM (\G 
)\rangle$ can be realized by division by $\LM (\G )$. So from a computational 
viewpoint, the effectiveness of Theorem 7.4 is done (the reader is referred to 
[G-I] for some examples of monomial algebras that satisfy the required conditions 
and for the algorithms written in pseudo-code). {\parindent=0pt\v5              
{\bf Finiteness of Global dimension}\par                                         
Let $\Omega\subset \B$ be a finite reduced subset. Following [An1] and [An2], 
Ufnarovski constructed in [Uf2] the {\it graph of chains} of $\Omega$ as a directed 
graph $\Gamma_{\rm C}(\Omega )$, in which the set of vertices $V$ is defined as
$$V=\{ 1\}\cup X\cup\{\hbox{all proper suffixes of}~u\in\Omega\},$$
and the set of edges $E$ contains all edges
$$1~\mapright{}{}~x_i~\hbox{for every}~x_i\in X$$
and edges defined by the rule
$$\begin{array}{rcl} u,v\in V-\{ 1\} ,~u~\mapright{}{}~v~\hbox{in}~E&\Leftrightarrow &
\hbox{there is a {\it unique}}~w\in\Omega ~\hbox{such that}\\
&{~~~~}&uv=\left\{\begin{array}{l} w,~\hbox{or}\\
sw,~s\in\B\end{array}\right.\end{array}$$ For $n\ge -1$, an $n$-{\it chain} in 
$\Gamma_{\rm C}(\Omega )$ is a word which can be read in the graph during a path of 
length $n+1$, starting from $1$. The set of all $n$-chains in $\Gamma_{\rm 
C}(\Omega )$ is denoted by $C_n$. For example, $C_{-1}=\{ 1\}$, $C_0=X$, and 
$C_1=\Omega$.\v5                                                                
{\bf 7.4. Theorem}  Suppose that $\G$ is a finite reduced Gr\"obner basis of $I$ 
with respect to $(\B ,\prec_{gr} )$. The following statements hold.\par          
(i) If the chain graph $\Gamma_{\rm C}(\LM (\G ))$ has no cycles, then 
gl.dim$(R/I)\le $ gl.dim$(R/\langle\HT (I)\rangle )\le $ gl.dim$(R/\langle\LM (\G 
)\rangle )<\infty$.\par                                                         
(ii) If the chain graph $\Gamma_{\rm C}(\LM (\G ))$ does not contain any 
$d$-chains, then gl.dim$(R/I)\le $ gl.dim$(R/\langle\HT (I)\rangle )\le $ 
gl.dim$(R/\langle\LM (I)\rangle )\le d$. \vskip 6pt                             
{\bf Proof} This follows from ([Uf2] Theorem 12), ([G-I2] Remark 3.14), ([An1] 
Theorem 4), the foregoing Theorem 4.8(ii) and section 5.\QED}\v5            
Consider any Gro\"bner basis $\G$ in the free $K$-algebra $K\langle 
X_1,X_2,X_3\rangle$ with $\LM (\G )=\{ X_2X_1, X_3X_1,X_3X_2\}$, then, the chain 
graph $\Gamma_{\rm C}(\LM (\G ))$ looks like
$$
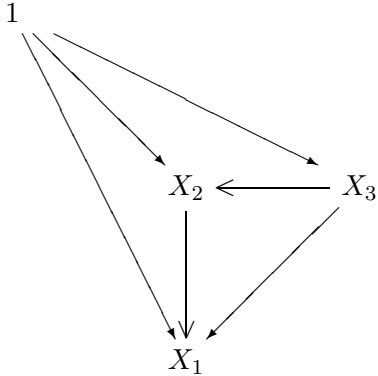
\begin{diagram}
1&                    &                      &                &\\
&\SE\ESE\SSE          &                      &                &\\
&                     &X_2                   &\lTo            &X_3\\
&                     &\dTo                  &\SW             &\\
&                     &X_1                   &                &\\
\end{diagram}$$
Clearly, $\Gamma_{\rm C}(\LM (\G ))$ does not contain 3-chains. Hence 
gl.dim$(R/\langle\G\rangle )\le $ gl.dim$(R/\langle\HT (\G )\rangle )\le$ 
gl.dim$(R/\langle\LM (\G )\rangle )\le 3$. In general, it follows from Theorem 7.4 
that the following result holds.{\parindent=0pt\v5                              
{\bf 7.5. Theorem}   Suppose that $\G$ is a finite reduced Gr\"obner basis of $I$ 
with respect to $(\B ,\prec_{gr} )$ such that
$$\LM (\G )=\left\{ X_jX_i~\Big |~1\le i<j\le n\right\} ,$$
(hence $R/I$ has a PBW $K$-basis). Then gl.dim$(R/I)\le $ gl.dim$(R/\langle\HT 
(I)\rangle )\le $ gl.dim$(R/\langle\LM (I)\rangle )\le n$.
\par\QED}\v5                                                
For the $K$-algebras  $A=R/I$ and $R/\langle\HT (I)\rangle$, to use Theorem 4.8(ii) 
and the results of section 5 effectively in examining whether they have finite 
global dimension, one can also check algorithmically if $K$ has a finite projective 
resolution over the $\mathbb{N}$-graded monomial algebra $\OV A=R/\langle\LM 
(I)\rangle=R/\langle\LM (\G )\rangle$ (for instance, the Anick resolution). 
Nowadays some well-developed computer algebra systems such as BERGMAN [CU] can 
produce a reduced Gr\"obner basis $\G$ for $I$ and  an Anick resolution for $K$. 
{\parindent=0pt\v5 {\bf Remark}                                                  
(i) In [G-IL], it was pointed out that V. Borisenko  proved in 1985 that a finitely 
presented monomial algebra $A=R/\langle\Omega\rangle$ satisfies some polynomial 
identity if and only if the Ufnarovski graph $\Gamma (\Omega )$ of $\Omega$ has no 
multiple vertex, and an algorithm was given to recognize this criterion. We have 
not yet explored what will happen to an algebra if its associated monomial algebra 
is a PI algebra.
\par
(ii)  We have not yet explored, on the basis of section 5, what will happen to an
algebra if its associated monomial algebra is a Koszul algebra. In [Li2], only very
little about this topic was discussed.\par                                     
(iii) We have not yet discussed any realization of the lifting properties for 
quotient algebras of commutative polynomial algebras and the coordinate rings of 
quantum affine $K$-spaces, either.} \par                                          
We leave these tasks for the successive work.

\v5 \centerline{References}{\parindent=1.2truecm\par
\re{[An1]} D. J. Anick, On Monomial algebras of finite global dimension, {\it
Trans. Amer. Math. Soc}., 1(291)(1985), 291--310.
\re{[An2]} D. J. Anick, On the homology of associative algebras, {\it Trans. Amer.
Math. Soc}., 2(296)(1986), 641--659.
\re{[Bu]} B. Buchberger, Gr\"obner bases: An algorithmic method in polynomial ideal
theory, in: {\it Multidimensional Systems Theory}, N.K. Bose, ed., Reidel,
Dordrecht, 1985, 184--232.
\re{[CU]} S. Cojocaru and V. Ufnarofski, BERGMAN under MS-DOS and Anick's
resolution, {\it Discrete Mathematics and Theoretical Computer Science}, 1(1997),
139--147.
\re{[G-I1]} T.~Gateva-Ivanova, Algorithmic determination of the Jacobson radical of
monomial algebras, in: {\it Proc. EUROCAL'85}, LNCS Vol. 378, Springer-Verlag,
1989, 355--364.                                         \re{[G-I2]}
T.~Gateva-Ivanova, Global dimension of associative algebras,
 in: {\it Proc. AAECC-6}, LNCS, Vol.357, Springer-Verlag, 1989, 213--229.
\re{[G-IL]} T.~Gateva-Ivanova and V.~Latyshev, On recognizable properties of
associative algebras, in: {\it Computational Aspects of Commutative Algebra}, from
a special issue of the Journal of Symbolic Computation, L.~Robbiano ed., Academic
Press, 1989, 237--254.
\re{[Gr]}  E.L. Green, {\it Noncommutative Gr\"obner Bases, A Computational \&
Theoretical Tool}, Lecture notes, Dec. 15, 1996,
www.math.unl.edu/~shermiller2/hs/green2.ps
\re{[GZ]} E. L. Green and D. Zacharia, The cohomology ring of a monomial algebra,
{\it Manuscripta Math}., 1(85)(1994), 11 -- 23.
\re{[K-RW]} A.~Kandri-Rody and V.~Weispfenning, Non-commutative Gr\"obner bases in
algebras of solvable type, {\it J. Symbolic Comput.}, 9(1990), 1--26.
\re{[Li1]} Huishi Li, {\it Noncommutative Gr\"obner Bases and Filtered-Graded
Transfer}, LNM, 1795, Springer-Verlag, 2002.
\re{[Li2]} Huishi Li, The general PBW property, {\it Alg. Colloquium}., in press.
arXiv:math.RA/0609172.
\re{[Li3]} Huishi Li, {\it Algebras and their Associated Monomial Algebras, ---
With Applications of Gr\"obner Bases}, Monograph, in preparation.
\re{[Mor]} T.~Mora, An introduction to commutative and noncommutative Gr\"obner
bases, {\it Theoretical Computer Science}, 134(1994), 131--173.
\re{[Nor]} P. Nordbeck, On the finiteness of Gr\"obner bases computation in
quotents of the free algebra, {\it Comm. and Comp}., 11(3)(2001), 157--180.
\re{[NVO]} C.~N$\check{\rm a}$st$\check{\rm a}$sescu and F.~Van Oystaeyen, {\it
Graded ring theoey}, Math. Library 28, North Holland, Amsterdam, 1982.
\re{[Uf1]} V. Ufnarovskii, A growth criterion for graphs and algebras defined by
words, {\it Mat. Zametki}, 31(1982), 465--472 (in Russian); English translation:
{\it Math. Notes}, 37(1982), 238--241.
\re{[Uf2]} V. Ufnarovskii, On the use of graphs for computing a basis, growth and
Hilbert series of associative algebras, (in Russian 1989), {\it Math. USSR
Sbornik}, 11(180)(1989), 417-428.

\end{document}